\documentclass[preprint,12pt]{elsarticle}

\usepackage{amssymb}
\usepackage{amsmath}
 \usepackage{amsmath,amssymb,amsfonts,amsthm}
 \usepackage{graphicx}
\usepackage{layout}
\usepackage[latin1]{inputenc}
\usepackage[french,english]{babel}
\usepackage[french]{babel}
\usepackage{venndiagram}
\usepackage{cancel}
\newtheorem{theorem}{Theorem}
\newtheorem{definition}{Definition}

\newtheorem{corollary}{Corollary}
\newtheorem{proposition}{Proposition}
\newtheorem{lemma}{Lemma}

\newtheorem{example}{Example}
\newtheorem{remark}{Remark}

\journal{Indagationes Mathematicae}

\begin{document}

\begin{frontmatter}


\title{Besselian Schauder Frames and the Structure of Banach Spaces}
 \author[label1]{Rafik Karkri}
 \affiliation[label1]{organization={FSK, Ibn Tofail University},
             addressline={P.O.Box 133},
             city={Kenitra},
             postcode={14000},
             country={Morocco}}
\ead{karkri.rafik@gmail.com}

\begin{abstract}
Schauder bases are fundamental tools for analyzing the structure of Banach spaces. In this work, we show that Besselian Schauder frames (BSF) play a similar role in certain contexts.  BSF are a new class of Schauder frames, lying between unconditional and general frames. We first prove that every unconditional Schauder frame (USF) is BSF, but the reverse implication is false. Specifically, we extend several well-known results of Karlin and James to Banach spaces with BSF, particularly to those with USF. We prove that many classical Banach spaces do not admit BSF, and in particular, do not admit USF.

Before establishing these results, for every Banach space $E$ with a finite dimensional decomposition, we provide an explicit method to construct a Schauder frame for $E$. In particular, Szarek's Banach space has a Schauder frame, which famously lacks a Schauder basis. This finding provides strong motivation for extending classical Schauder basis theory to the framework of Schauder frames.
\end{abstract} 
\begin{keyword}
Schauder frame \sep frame \sep weakly sequentially complete Banach space \sep Besselian Schauder frame \sep unconditional Schauder frame \sep unconditional Schauder basis.
\MSC[2008] 46B04 \sep 46B10 \sep 46B15 \sep 46B25 \sep 46B45 \sep
46A32
\end{keyword}
\end{frontmatter}



\section{Introduction}
In 1932, Banach and Mazur \cite[theorem 9, p. 185]{S.Banach.1932} 
proved that every separable Banach space is isometrically 
embeddable into $C([0,1])$ (the Banach space of continuous functions on the closed
interval $ [0,1]$).
The existence of a Schauder basis for $C([0,1])$ 
\cite[example 4.1.11, p. 352]{R.Megginson.1998} led Banach to ask whether every separable Banach space has a Schauder basis \cite[p. 111]{S.Banach.1932}.

Recall that a sequence
$\left( x_{n} \right) _{n\in \mathbb{N}^{*}}$ in a Banach space 
$E$ is called a Schauder basis (in short, basis) for $E$ if, for every $x\in E$, there exists a unique sequence of scalars $\left( \alpha_{n} \right) _{n\in \mathbb{N}^{*}}$
such that 
$x=\sum_{n=1}^{\infty}\alpha_{n}x_{n}$. If, in addition, the series
 $\sum_{n}\alpha_{n}x_{n}$  converges unconditionally to $x$ for every 
$x\in E$, the sequence is called an unconditional Schauder basis (in short, unconditional basis) for $E$.

This problem (known as the basis problem) remained open for a long time and was solved in 1973 in the
negative by P. Enflo \cite{P.Enflo.1973}.

In 1987, S. J. Szarek \cite{S.J.Szarek} constructed a complemented subspace without a Schauder basis in a space that does possess a Schauder basis, thereby solving a fifty-year-old problem. This result established that the Schauder basis property is not inherited by complemented subspaces.

In 1948, Karlin  \cite{S.Karlin.1948} 
used the concept of unconditional bases to investigate the relationships between the structures of a Banach space $E$, its dual $E^{*}$ and its bidual $E^{**}$.
Among other results, he studied reflexivity and weak sequential completeness, proving the following:
\begin{quote}
\textbf{Theorem (Karlin):}
Let $E$ be a Banach space with an unconditional basis. If $E^{*}$ is separable, then:
\begin{enumerate}
\item
$E^{*}$ is weakly sequentially complete;
\item
$E^{*}$ has an unconditional basis.
\end{enumerate}
\end{quote}
In the same paper, Karlin posed the following question:
\begin{quote}
\textbf{Question (Karlin):}
  If $E^{*}$ has a  basis, does $E$ itself have a basis?
  \end{quote}
  
 In 1971, Johnson, Rosenthal, and Zippin
 \cite{Johnson.Rosenthal.Zippin} provided an affirmative answer to this question, thereby resolving one of the fundamental problems in basis theory.

In $1950$, James  \cite{James.1950}  used the concepts of shrinking and boundedly complete
basis as tools in the structure theory of Banach spaces.

The concept of frames (generalizing the concept of bases) for Hilbert spaces was introduced in $1952$ by Duffin and Schaeffer   \cite{Duffin.1952}.  A  frame for a Hilbert space  
 $\left( \mathcal{H},\left \langle \cdot,\cdot\right \rangle \right) $ is a sequence
 $\left(h_{n}\right)_{n\in \mathbb{N}^{* }}$ in $\mathcal{H}$ for which there exist constants $A,B>0$ such that: 
\begin{equation}
\label{std.frame}
A\left \Vert h\right \Vert _{\mathcal{H}}^{2}\leq \underset{n=1}{\overset{\infty }{\sum }}
\left \vert \left \langle h,h_{n}\right \rangle \right \vert ^{2}
\leq B\left \Vert h\right \Vert _{\mathcal{H}}^{2},\; \;(h\in \mathcal{H}).
\end{equation}
This condition allows each    
$x\in\mathcal{H}$ to be 
written as:
 \[ x=\underset{n=1}{\overset{\infty}{\sum }}
 \left \langle x,h_{n}\right \rangle h_{n},\]
 but the  linear independence between the elements of 
 $\left(  h_{n}\right)_{n\in \mathbb{N}^{*}}$ is not required.
If the upper inequality in \eqref{std.frame} holds,
$\left(h_{n}\right)_{n\in \mathbb{N}^{*}}$ is called 
a Bessel sequence with bound $B$.
 
The concept of frames (as a generalization 
	  of Schauder bases and  Hilbert frames) for  Banach spaces was 
 introduced in $2008$ by 
 Casazza, Dilworth, Odell, Schlumprecht, and Zsak \cite{Casazza.D.S.Z.2008}.
A Schauder frame for a given  Banach space $E$ is  a sequence  
$\left((a_{n},b_{n}^{*})\right)_{n\in \mathbb{N}^{*}}\subset E\times E^{*}$
 such that, for each $x$ in $E$, we have:
 \[x=\underset{n=1}{\overset{\infty}{\sum }}
 b_{n}^{*}\left(x\right) a_{n}.\]

In 2021,  a separate group of authors \cite{karkri.zoubeir.2021} introduced the concept of Besselian Schauder frames for Banach spaces by adding a Besselian property to this definition. Specifically:
 \begin{quote}
\textbf{Definition.}
A sequence $\mathcal{F}:=\left(\left( a_{n},b_{n}^{*}\right)\right)_{n\in \mathbb{N}^{*}}$ in $E\times E^{*}$ is called a Besselian Schauder frame (BSF) for $E$ if the following conditions hold:
\begin{enumerate}
\item
$\mathcal{F}$ is a Schauder frame for $E$; 
\item
There exists a constant $B>0$ such that:
\begin{equation*}
\underset{n=1}{\overset{\infty}{\sum }}\left \vert b_{n}^{*}\left(x\right)\right \vert\left\vert y^{*}\left(a_{n}\right)\right \vert
\leq B\left\Vert x\right\Vert_{E}\left\Vert y^{*}\right\Vert
_{E^{*}},\;\;\;(x\in E,y^{*}\in E^{*}).
\end{equation*}
\end{enumerate} 
\end{quote}
The second condition generalizes the upper inequality in \eqref{std.frame}(but, in contrast to the Hilbert space case, it does not generally imply unconditional convergence of the frame) and reveals the dual space structure that is not apparent in the definition of frames for Hilbert spaces
(since $\mathcal{H}^{*}\cong \mathcal{H}$).

In 2009, R. Liu \cite{Rui.Liu.2009} extended James's results on shrinking and boundedly complete bases to unconditional Schauder frames by introducing the concepts of locally shrinking frames and locally boundedly complete frames.

In the same year, a separate group of authors---Carando and Lassalle \cite{Carando.Lassalle.2009}, and D. Carando , S. Lassalle and P. Schmidberg \cite{Carando.Lassalle.Schmidberg.2009}---extended these results to Banach spaces with unconditional atomic decompositions, and hence to Banach spaces with unconditional Schauder frames, by introducing a more general definition of a shrinking frame.
   
In this paper, we generalize these results (using a new and general definition of shrinking frame, which coincides with the natural definition in the basis setting) to the large class of Banach spaces with Besselian Schauder frames, and in particular to Banach spaces with unconditional Schauder frames. We prove the following theorem:
\begin{quote}
\textbf{Theorem A.}
Let $\mathcal{F}:=\left( \left( a_{n},b_{n}^{*}\right) \right)
_{n\in \mathbb{N}^{*}}$ be a Besselian Schauder frame for a Banach space $E$. Then:
\begin{enumerate}
\item
$\mathcal{F}$
is shrinking if and only if there is no subspace of $E$   isomorphic to 
$\ell_{1}$.
\item
$\mathcal{F}$ is boundedly complete   for $E$ 
if and only if 
$E$ is weakly sequentially complete
if and only if
there is no subspace of $E$ isomorphic to $c_{0}$.
\item
$E$ is reflexive 
 if and only if $E$ contains no isomorphic copies of $\ell_{1}$ or $c_{0}$.
\end{enumerate} 
\end{quote}
On the other hand, we generalize some well-known results of Karlin \cite{S.Karlin.1948}  to Banach spaces with Besselian Schauder frames (and, in particular, to those with unconditional Schauder frames); see \textbf{theorem B} in section \eqref{Notations.preliminaries}.

It is well-known that:
\begin{quote}
\begin{enumerate}
\item
Every Banach space with an unconditional Schauder basis possesses Pe\l czy\'{n}ski's property (u) \cite[proposition 3.5.3, p. 67]{Albiac.Kalton.2006}.
\item
The classical Banach spaces such as
 $C(\omega^{\omega})$, $C([0,1])$, 
$L_{1}([0,1])$, the James space $\mathcal{J}$,
$K(\ell_{p},\ell_{q})$,
$\ell_{p}\otimes_{\pi}\ell_{q^{*}}$ and $\ell_{p^{*}}\otimes_{\varepsilon}\ell_{q}$ (for $p\geq q>1$) do not admit unconditional Schauder bases.
\end{enumerate}
\end{quote}  
In this paper, we extend these results to Banach spaces with Besselian Schauder frames (and, in particular, to those with unconditional Schauder frames).

Finally, for every Banach space $E$ with a finite dimensional decomposition, 
 we provide an explicit method to construct a Schauder frame for $E$. As a result, Szarek's famous example \cite{S.J.Szarek} serves as an illustration of a Banach space with a Schauder frame that does not admit a Schauder basis. This example underscores the importance of generalizing classical results from Banach spaces with Schauder bases to those with Schauder frames.
\section{Notations and preliminaries}
\label{Notations.preliminaries}
Throughout this paper, $E$ denotes a Banach space (real or complex), and $E^{*}$ represents its topological dual.
\begin{enumerate}
\item
Definitions:
\begin{itemize}
\item 
Let $\mathcal{F}:=\left( \left( a_{n},b_{n}^{*}\right) \right)
_{n\in \mathbb{N}^{*}}$ be a fixed sequence in $E \times E^{*}$.
\item 
The dual sequence of $\mathcal{F}$ in $ E^{*} \times E^{**}$, is denoted by:
$$\mathcal{F}^{*}:=\left( \left( b_{n}^{*},J_{E}(a_{n})\right)\right) _{n\in \mathbb{N}^{*}},$$ 
 where $J_{E}:E\rightarrow E^{**}$ is the canonical embedding of $E$ into its bidual.
\item
The closed unit ball of $E$ is denoted by:
$$\mathbb{B}_{E}:=\left \{ x\in E:\left \Vert x\right \Vert _{E}\leq 1\right\}
.$$
\item
For $p \in ]1, +\infty[$, the conjugate exponent is:
$$p^{*}=\dfrac{p}{p-1}.$$
\item
Denote by $\mathcal{S}$ the set of all sign sequences:  
\[
\mathcal{S} = \left\{s=\left (s_{m}\right )_{m \in\mathbb{N}^{*}}\subset\mathbb{K}: \vert s_{m}\vert=1\right\}.
\]  
We equip $\mathcal{S}$ with the metric:  
\[
d(s,t)=\sum_{m=1}^\infty \frac{|s_m - t_m|}{2^m}.
\]  
For each component space $\{z \in \mathbb{K}: \vert z\vert= 1\}$ is compact (closed and bounded in $\mathbb{K}$). By Tychonoff's theorem, the product of countably many compact sets is compact in the product topology. The metric $d$ induces this topology since $s^{(k)} \to s$ in $(\mathcal{S},d)$ if and only if $s_{m}^{(k)} \to s_{m}$ in $\mathbb{K}$ for all $m$. Thus, $\mathcal{S}$ is compact.  
\end{itemize}
\item
Classical Banach spaces
\begin{itemize}
\item
We denote by $c_{0}$, $\ell_{p}$ (for $p\geq 1$) and $\ell_{\infty}$ the classical sequence spaces.
\item
 We denote by $L_{p} \left([0,1] \right)$ (for $p\geq 1$) the Banach space of equivalence classes of Lebesgue integrable functions 
 $f: [0,1]\rightarrow \mathbb{K}$, where $\mathbb{K}=\mathbb{R}$ or  $\mathbb{K}=\mathbb{C}$.
 \item
We denote by $C(\omega^{\omega})$ the Banach space of all continuous functions 
on the (compact metrizable) space $\omega^{\omega}$ of all ordinal numbers 
$\leq\omega^{\omega}$ endowed with the order topology.
 \item
We denote by $L\left( E\right) $ the Banach space of all bounded linear
operators $T: E\rightarrow E$.
\item
Let $F$ be a Banach space. We denote by $K(E,F)$ the Banach space of all compact operators $T: E\rightarrow F$.
\end{itemize}
\item
Tensor product of Banach spaces \cite{R.Schatten.1950}
\begin{itemize}
\item
We denote by
$E\otimes F:=span\left \{x\otimes y: (x,y)\in E\times F\right\}$
 the algebraic tensor product of $E$ and $F$, where $x\otimes y$ 
 is the bounded operator $E^{*}\rightarrow F$ defined by $f\mapsto f(x)y$. 
 \item
 A norm $\alpha$ on $E\otimes F$ is called cross norm if
 $\alpha (x\otimes y )=\Vert x \Vert_{E} \Vert y \Vert_{F}$ for all $x\in E$ and $y\in F$.
 \item
 A crossnorm $\alpha$ on $E\otimes F$ is called uniform if for all bounded linear operators $S\in L(E)$ and $T\in L(F)$, the linear operator: 
\begin{align*}
S\otimes T:\;  E\otimes F \longrightarrow E\otimes F ,\;\;
\sum_{r=1}^{R} x^{r}\otimes y^{r}\longmapsto\sum_{r=1}^{R}S(x^{r})\otimes T(y^{r})
\end{align*}
is bounded with respect to $\alpha$ and its norm satisfies:
\begin{align*}
\left\Vert S\otimes T\right\Vert_{L(E\otimes F)}
     &\leq 
     \Vert S\Vert _{L(F)}\Vert T\Vert_{L(E)}.
\end{align*}
Equivalently, for every tensor $u\in E\otimes F$,
\begin{align*}
\alpha\left(S\otimes T(u)\right)
     &\leq 
     \Vert S\Vert _{L(F)}\Vert T\Vert_{L(E)}
 \alpha\left(u\right).
\end{align*}
 \item 
 We denote by $  E\otimes_{\alpha}F$ the completion of $E\otimes F$ with 
 respect to the cross norm $\alpha $.
 \item 
The space $E\otimes_{\varepsilon}F$ is called  the injective  tensor product
 of $E$ and $F$, where the cross norm $\varepsilon$ is defined by:
 \[ \varepsilon\left (\underset{i=1}{\overset{n}{\sum }} x_{i}\otimes y_{i}
 \right )
 =\underset{\left(f^{*},g^{*}\right)
 \in \mathbb{B}_{E^{*}}\times \mathbb{B}_{F^{*}}}{\sup }
 \left \vert \underset{i=1}{\overset{n}{\sum }}
   f^{*}(x_{i})g^{*}(y_{i})\right  \vert.
  \]
\item
The space $E\otimes_{\pi}F$ is called  the projective tensor product
 of $E$ and $F$, where the cross norm $\pi$ is defined by:
\[ \pi(u)=inf\left \{\underset{i=1}{\overset{ n}{\sum }}\Vert x_{i}\Vert_{E}    \Vert y_{i}\Vert_{F}:
u=\underset{i=1}{\overset{n}{\sum }} x_{i}\otimes y_{i} \right\},\]
the infimum being taken over all the representations of $u$ in $  E\otimes F$.
\end{itemize}
\item
Weak convergence properties in Banach spaces
\begin{itemize}
\item \cite[definition 2.3.9, p. 38]{Albiac.Kalton.2006}.
A sequence $\left(x_{n}\right)_{n\in \mathbb{N}^{*}}$ in $E$ is called weakly
 Cauchy if $\lim_{n\rightarrow\infty}y^{*}(x_{n})$ exists for every $y^{*}\in E^{*}$.
\item \cite[definition 2.4.3, p. 39]{Albiac.Kalton.2006}.
A series $\sum_{n} x_{n}$ in $E$ is called weakly
unconditionally Cauchy (in short, WUC) if $\sum_{n}
\left \vert y^{*}\left(x_{n}\right) \right \vert $ converges for
every $y^{*}\in E^{*}$.
\item \cite[definition 3.5.1, p. 66]{Albiac.Kalton.2006}.
$E$ has property (u) if for every weakly Cauchy sequence $\left(x_{n}\right) _{n\in \mathbb{N}^{*}}$ in $E$, there exists a WUC series $\sum_{k}u_{k}$ in
$E$ such that $x_{n}-\sum_{k=1}^{n}u_{k}\rightarrow 0$ weakly.
\item  
\cite[p. 38]{Albiac.Kalton.2006}. 
 $E$ is called weakly sequentially complete
(in short,  WSC) if every weakly Cauchy sequence in $E$ converges weakly to some limit in $E$.
\end{itemize}
\item
Approximation properties in Banach spaces
\begin{itemize}
\item
\cite[definition 3.4.26, p. 330]{R.Megginson.1998}.
We say that $E$ has the approximation property (in short, AP) if,
for every Banach space $F$, the set of finite-rank operators is
dense in $K(F,E)$.
\item
\cite[theorem 3.4.32, p. 334]{R.Megginson.1998}. (A. Grothendieck, 1955).
We say that $E$ has the approximation property if 
for every compact set $K\subset E$ and every $\varepsilon>0$,
there exists a finite-rank operator $T$ of $L(E)$ such that:
$$\left \Vert T(x)-x \right\Vert_{E}< \varepsilon , \;\;(x\in K).$$  
\item
\cite[definition 4.1.34, p. 365]{R.Megginson.1998}.
We say that  $E$ has the bounded approximation property 
(in short, BAP) if there exists a constant $t>0$ having the property that, for every compact set $K\subset E$ and every $\varepsilon>0$,
 there exists a finite-rank operator $T$ of $L(E)$ such that:
$$\left \Vert T\right \Vert_{L(E)}\leq  t\;\; \text{and}\;\;
\left \Vert T(x)-x \right \Vert_{E}< \varepsilon ,\;\;(x\in K).$$ 
\end{itemize}
\item 
Finite dimensional decomposition \cite{Johnson.Rosenthal.Zippin}.\\
$E$ is said to have a finite dimensional decomposition
(in short,  FDD) if there exists a sequence $(E_{n})_{n\in\mathbb{N}^{*}}$ of finite dimensional subspaces of $E$ such that every $x\in E$ admits a unique representation:
\[
x = \sum_{n=1}^\infty P_n(x), \quad \text{where } P_n(x) \in E_n \text{ for all } n \in \mathbb{N}^*.
\]
In this case, each $P_n \colon E \to E_n$ is a bounded linear projection.
\end{enumerate}
To simplify the summary of the main findings of this paper,  we use the following notations.\\
\textbf{Notation 1.}
\begin{enumerate}
\item $\mathcal{F}(SF)$: $\mathcal{F}$ is a Schauder frame.
\item $\mathcal{F}(USF)$: $\mathcal{F}$ is an unconditional Schauder frame.
\item $\mathcal{F}(BSF)$: $\mathcal{F}$ is a Besselian Schauder frame.
\item $\mathcal{F}(SB)$: $\mathcal{F}$ is a Schauder basis.
\item $\mathcal{F}(USB)$: $\mathcal{F}$ is an unconditional  Schauder basis.
\item $\mathcal{F}(BSB)$: $\mathcal{F}$ is a Besselian Schauder basis.
\end{enumerate}
In this paper, we establish the following relationships. (Remark: All the relationships below hold in the bases setting as well, except possibly item \eqref{dual.bess.WSC.item}.)\\
\textbf{Theorem B.}
\begin{enumerate}
\item
$\mathcal{F}(USB)\Longleftrightarrow \mathcal{F}(BSB)\Longrightarrow
\mathcal{F}(USF) \Longrightarrow \mathcal{F}(BSF)\Longrightarrow\mathcal{F}(SF).
$
\item
$\mathcal{F}(SF) \cancel{\Longrightarrow} \mathcal{F}(BSF)
\cancel{\Longrightarrow} \mathcal{F}(USF)
\cancel{\Longrightarrow} \mathcal{F}(USB)$
\item
$\mathcal{F}(SB) \Longrightarrow \mathcal{F}(SF)
$ and 
$\mathcal{F}(SF) \cancel{\Longrightarrow} \mathcal{F}(SB)$.
\item
If $E$ is WSC, then $\mathcal{F}(BSF)\Longleftrightarrow \mathcal{F}(USF)$.
\item
$E$ has a $\mathcal{F}(SF)$ $ \cancel{\Longrightarrow} $  $E^{*}$ has a $\mathcal{F}(SF)$.
\item
$E$ has a $\mathcal{F}(BSF)$ and $E^{*}$ is separable
$\Longrightarrow$ $E^{*}$ is WSC and has an $\mathcal{F}(USF)$.
\item
\label{dual.bess.WSC.item}
$E^{*}$ has a $\mathcal{F}(BSF)$
$ \Longrightarrow$  $E^{*}$ is WSC.
\item
$E$ and $F$ have  $\mathcal{F}(BSF)$ $ \cancel{\Longrightarrow} $ $  E\otimes_{\pi}F$ has a $\mathcal{F}(BSF)$.
\item
$E$ and $F$ have  $\mathcal{F}(BSF)$ $ \cancel{\Longrightarrow} $ $  E\otimes_{\varepsilon}F$ has a $\mathcal{F}(BSF)$.
\item
$E^{*}$ has a $\mathcal{F}(SF)$ $ \Longrightarrow $ $E$ has a $\mathcal{F}(SF)$.
\item
$E^{*}$ has a $\mathcal{F}(BSF)$ $ \cancel{\Longrightarrow}$ $E$ has a $\mathcal{F}(BSF)$.
\item
$E^{*}$ has an $\mathcal{F}(USF)$ $ \cancel{\Longrightarrow}$ $E$ has an $\mathcal{F}(USF)$.
\end{enumerate}
\section{Schauder frames}
The purpose of this section is to examine the relationships between the following concepts:
\begin{itemize}
\item
Schauder frames (and bases),
\item
Besselian Schauder frames (and bases), and
\item
unconditional Schauder frames (and bases).
\end{itemize}
Additionally, we establish several characterizations of Besselian Schauder frames.
\begin{definition} \cite[definition 2.2, p. 68]{Casazza.D.S.Z.2008}
A sequence $\mathcal{F}$ is called a Schauder frame (in short, frame) (resp.
unconditional Schauder frame (in short, unconditional frame) ) for $E$ if for all $x\in E$, the series $\sum_{n}
b_{n}^{*}\left( x\right) a_{n}$ is convergent (resp. unconditionally
convergent) in $E$ to $x$.\\

It is an open problem whether every complemented subspace of a space with
unconditional basis has an unconditional basis
 \cite[problem 1.8, p. 279]{Casazza.AP.2001}.
 The example \eqref{frame.complemented} shows (without any difficulty)
 that a frame passes to complemented subspaces.
\begin{example}
\label{frame.complemented}
Every complemented subspace of a Banach space with a frame (resp. an unconditional frame) has a frame (resp. an unconditional frame).
\end{example}
\textbf{Proof.}
Let $E$ be a Banach space with an unconditional frame
$\mathcal{F}:=\left( \left(a_{n},b_{n}^{*}\right)
\right) _{n\in \mathbb{N}^{*}}\subset E\times E^{*}$,
and let $F$ be a complemented subspace of $E$ with corresponding continuous projection  $P:E\rightarrow F$. For any permutation $\sigma $ of $\mathbb{N}^{*}$ and any $x\in F$, we have:
$$
x =P\left( x\right) 
=\overset{\infty}{\underset{n=1}{\sum }}
b_{\sigma(n)}^{*}(x)P\left(a_{\sigma(n)}\right).
$$
It follows that 
$\left( \left( P\left(a_{n}\right),
b_{n\left \vert F\right. }^{*}\right) \right)_{n\in \mathbb{N}^{*}}$
 is an unconditional  frame for $F$ 
 (where $b_{n\left \vert F\right. }^{*}$
 is the restriction of $b_{n}^{*}$ to $F$). 
 If $\mathcal{F}$ is just a frame, it suffices to put ($\sigma=Id$). \qed
\end{definition}
In \cite[proposition 6.5, p. 300]{Casazza.AP.2001}, it is proved  that if $(E_{n})_{n\in\mathbb{N}^{*}}$ is a FDD of $E$ and each 
$E_{n}$ has a basis $((e_{n,i},e^{*}_{n,i}))_{i=1}^{d_{n}}$ with basis constant $K_{n}$
such that $\underset{n}{\sup}K_{n}<\infty$ then, the sequence
$(((e_{n,i},e^{*}_{n,i}))_{i=1}^{d_{n}})_{n\in\mathbb{N}^{*}}$ 
is a  basis for $E$. 

In the  frames setting (proposition \eqref{FDD.impl.frame}),
we provide an explicit method to construct a  frame for $E$ from the
bases of each $E_{n}$ without requiring
$\underset{n}{\sup}K_{n}<\infty$.
\begin{proposition}
\label{FDD.impl.frame}
Suppose that $E$ has an FDD. Then $E$ has a  frame.
\end{proposition}
\textbf{Proof.}
Let $(E_{n})_{n\in\mathbb{N}^{*}}$ be an FDD of $E$. Suppose that  for each
$n\in\mathbb{N}^{*}$, the sequence 
$\left((e_{n,i},e^{*}_{n,i})\right)_{i=1}^{d_{n}}$ is a basis for  $E_{n}$ with
$\Vert e^{*}_{n,i} \Vert=\Vert e_{n,i} \Vert=1$.
Define a sequence $((a_{n,j}, b^{*}_{n,j}))_{j=1}^{d_{n}^{2}}$ in $E_{n}$ by:
 $$a_{n,kd_{n}+r}=d_{n}^{-1}e_{n,r}\;\;\text{and}\;\;b^{*}_{n,kd_{n}+r}=e^{*}_{n,r},$$
 for each 
 $1\leq r \leq d_{n}$ and $0 \leq k \leq d_{n}-1$. Then 
 $((a_{n,j}, b^{*}_{n,j}))_{j=1}^{d_{n}^{2}}$   is a  frame for 
$E_{n}$. We claim that the sequence 
 $\left(\left((a_{n,j},b^{*}_{n,j}\circ P_{n})\right)_{j=1}^{d_{n}^{2}}\right)_{n\in\mathbb{N}^{*}}$
 is a  frame for $E$. 
To verify this, let 
  $x\in E$ and $\varepsilon >0$. By the properties of the FDD, there exists $N\in\mathbb{N}^{*}$ such that for all $n\geq N$, we have
  $$\left\Vert \overset{n-1}{\underset{k=1}{\sum }}P_{k}(x)-x\right\Vert_{E}\leq\dfrac{\varepsilon}{2}\;\;\;\text{and}\;\;\;
  \left\Vert P_{n}(x)\right\Vert_{E}\leq \dfrac{\varepsilon}{2}.$$
Hence it suffices to prove that for each $1\leq m < d_{n}^{2}$ we have
 $$\left\Vert S_{n-1,m}(x)-x\right\Vert_{E}\leq \varepsilon ,$$ 
 where
 $$S_{n-1,m}(x):= \overset{n-1}{\underset{k=1}{\sum }}
 \overset{d_{k}^{2}}{\underset{j=1}{\sum }}
 b^{*}_{k,j}(P_{k}(x))a_{k,j}+\overset{m}{\underset{i=1}{\sum }}
 b^{*}_{n,i}(P_{n}(x))a_{n,i}.$$
For each $1\leq m < d_{n}^{2}$, there exists $0\leq\alpha \leq d_{n}-1 $ such that:
$$\alpha d_{n} < m \leq (\alpha +1)d_{n}.$$ 
We then decompose $A:=\left \Vert\overset{m}{\underset{i=1}{\sum }}
 b^{*}_{n,i}(P_{n}(x))a_{n,i}\right \Vert_{E}$ as follows:
\begin{align*}
 A
 &=
 \left \Vert\overset{\alpha d_{n}}{\underset{i=1}{\sum }}
 b^{*}_{n,i}(P_{n}(x))a_{n,i}
 +\overset{m}{\underset{i=\alpha d_{n}+1}{\sum }}
 b^{*}_{n,i}(P_{n}(x))a_{n,i}\right \Vert_{E}\\
 &=
 \left \Vert \sum_{k=0}^{d_{n}-1}\sum_{r=1}^{d_{n}}
 b^{*}_{n,kd_{n}+r}(P_{n}(x))a_{n,kd_{n}+r}
 +\overset{m-\alpha d_{n}}{\underset{r=1}{\sum }}
 b^{*}_{n,\alpha d_{n}+r}(P_{n}(x))a_{n,\alpha d_{n}+r}\right \Vert_{E}\\
 &=
 \left \Vert \sum_{k=0}^{d_{n}-1}\sum_{r=1}^{d_{n}}
 d_{n}^{-1}e^{*}_{n,r}(P_{n}(x))e_{n,r}
 +\overset{m-\alpha d_{n}}{\underset{r=1}{\sum }}
 d_{n}^{-1}e^{*}_{n,r}(P_{n}(x))e_{n,r}\right \Vert_{E}\\
 &=
 \left \Vert \sum_{k=0}^{d_{n}-1}
 d_{n}^{-1}P_{n}(x)
 +\overset{m-\alpha d_{n}}{\underset{r=1}{\sum }}
 d_{n}^{-1}e^{*}_{n,r}(P_{n}(x))e_{n,r}\right \Vert_{E}\\
 &\leq 
 \left \Vert P_{n}(x)\right \Vert_{E}
 +\frac{m-\alpha d_{n}}{d_{n}}\left \Vert P_{n}(x)\right \Vert_{E}.
\end{align*}
Consequently
$$
\left\Vert S_{n-1,m}(x)-x\right\Vert_{E} 
 \leq  \left\Vert \overset{n-1}{\underset{k=1}{\sum }}
 P_{k}(x)
 -x\right\Vert_{E}+ \left\Vert\overset{m}{\underset{i=1}{\sum }}
 b^{*}_{n,i}(P_{n}(x))a_{n,i}\right\Vert_{E}\leq\varepsilon,
$$
which finishes the proof.\qed
\begin{example}
\label{example.frame.FDD}
Let  $\left( \left(e_{n},e_{n}^{*}\right) \right) _{n\in \mathbb{N}^{*}}$ 
be the canonical basis  for  
$c_{0}$. Since 
$(span(e_{n}))_{n\in\mathbb{N}^{*}}$ is an FDD of $c_{0} $,
proposition \eqref{FDD.impl.frame} implies that the sequence 
$$\left(
\left(\frac{1}{2}e_{1},e_{1}^{*}\right), \left(\frac{1}{2}e_{1},e_{1}^{*}\right), \left(\frac{1}{2}e_{2},e_{2}^{*}\right), 
\left(\frac{1}{2}e_{2},e_{2}^{*}\right),
\left(\frac{1}{2}e_{3},e_{3}^{*}\right),
\left(\frac{1}{2}e_{3},e_{3}^{*}\right),...\right )$$
is a frame (which is not a basis) for
$c_{0}$. 
\end{example}
\begin{example}
\label{example.uncond.frame.FDD}
The sequence 
$$\left(
\left(\frac{1}{2}e_{1},e_{1}^{*}\right), \left(\frac{1}{2}e_{1},e_{1}^{*}\right), \left(\frac{1}{2}e_{2},e_{2}^{*}\right), 
\left(\frac{1}{2}e_{2},e_{2}^{*}\right),
\left(\frac{1}{2}e_{3},e_{3}^{*}\right),
\left(\frac{1}{2}e_{3},e_{3}^{*}\right),...\right )$$
is an unconditional frame for $c_{0}$. 
\end{example}
\textbf{Proof.}
Let $x\in c_{0}$. It suffices to apply the classical characterization of
  the unconditionally convergent series  \cite[lemme 1.45 (4), p. 31]{A.Abramovich.2002} to the series $\sum_{n}2^{-1}e_{n}^{*}(x)e_{n}$ (which converges to $2^{-1}x$).\qed
\begin{corollary}
There exists a separable Banach space without a basis which has a frame.
\end{corollary}
\textbf{Proof.}
By a result of S.J. Szarek \cite{S.J.Szarek}, there is 
a Banach spaces $E$ with an FDD  which does not have  a  basis.  Then  by
Proposition \eqref{FDD.impl.frame} this space has a frame.\qed
 \begin{definition}\cite{karkri.zoubeir.2021}
 \label{bess.sequence}
A sequence $\mathcal{F}$ is called Besselian  if
there is a constant $B>0$ such that: 
\begin{equation*}
\underset{n=1}{\overset{\infty}{\sum }}\left \vert b_{n}^{*}\left(
x\right) \right \vert \left \vert y^{*}\left( a_{n}\right) \right \vert
\leq B\left \Vert x\right \Vert _{E}\left \Vert y^{*}\right \Vert
_{E^{*}},\;\;\;(x\in E,y^{*}\in E^{*}).
\end{equation*}
 \end{definition}
 \begin{remark}
For a Besselian sequence $\mathcal{F}$ of $E$, the quantity 
\begin{equation*}
\mathcal{L}_{\mathcal{F}}:=\underset{\left(x,y^{*}\right)
 \in \mathbb{B}_{E}\times \mathbb{B}_{E^{*}}}{\sup }
 \left( \underset{n=1}{\overset{\infty }{\sum }}
 \left \vert b_{n}^{*}\left(x\right) \right \vert \left
\vert y^{*}\left(a_{n}\right) \right \vert \right)
\end{equation*}
is finite and for each $ x\in E$ and $y^{*}\in E^{*}$, 
the following inequality holds
\begin{equation*}
\underset{n=1}{\overset{\infty }{\sum }}\left \vert b_{n}^{*}\left(
x\right) \right \vert \left \vert y^{*}\left( a_{n}\right) \right \vert
\leq \mathcal{L}_{\mathcal{F}}\left \Vert x\right \Vert _{E}\left \Vert
y^{*}\right \Vert_{E^{*}}.
\end{equation*}  
\end{remark}
\begin{definition}\cite{karkri.zoubeir.2021}
A sequence $\mathcal{F}$ is called a Besselian Schauder frame
(in short, Besselian frame) for $E$ if it is simultaneously:
\begin{itemize}
\item
A Besselian sequence (satisfying definition \eqref{bess.sequence}), and
\item
A Schauder frame for $E$.
\end{itemize}
\end{definition}
\begin{example}
\label{example.bess.frame.c_0}
The canonical basis $\left( \left(
e_{n},e_{n}^{* }\right) \right) _{n\in \mathbb{N}^{*}}$ 
of $c_{0}  $ is a Besselian frame for  
$c_{0} $.
\end{example}
\textbf{Proof.}
Let \(y^* \in c_0^*\) and \(x = \sum_{n=1}^\infty x_n e_n \in c_0\). We identify \(c_0^*\) with \(\ell_1\) via the isometric isomorphism \(y^* \mapsto (y_n)_{n=1}^\infty\), where \(y_n = y^*(e_n)\) and \(\|y^*\|_{c_0^*} = \sum_{n=1}^\infty |y_n|\). Then \(e_n^*(x) = x_n\), and we have:

\[
\sum_{n=1}^\infty |e_n^*(x)| \, |y^*(e_n)| = \sum_{n=1}^\infty |x_n| \, |y_n|.
\]
Since \(|x_n| \leq \|x\|_{c_0}\) for all \(n\), it follows that:
\[
\sum_{n=1}^\infty |x_n| \, |y_n| \leq \|x\|_{c_0} \sum_{n=1}^\infty |y_n| = \|x\|_{c_0} \, \|y^*\|_{c_0^*}.
\]
Thus $\left( \left(
e_{n},e_{n}^{* }\right) \right) _{n\in \mathbb{N}^{*}}$ is a Besselian frame for \(c_0\).\qed
\begin{example}
\label{example.bess.frame}
The canonical basis 
$\left( \left(e_{n},e_{n}^{* }\right) \right) _{n\in \mathbb{N}^{*}}$ 
of $\ell_{1}  $ is a Besselian frame for  
$\ell_{1} $.
\end{example}
\textbf{Proof.}
\label{bsf.exmample.l1}
Let \(y^* \in \ell_1^*\) and \(x = \sum_{n=1}^\infty x_n e_n \in \ell_1\). We identify \(\ell_1^*\) with \(\ell_\infty\) via the isometric isomorphism \(y^* \mapsto (y_n)_{n=1}^\infty\), where \(y_n = y^*(e_n)\) and \(\|y^*\|_{\ell_1^*} = \sup_{n \in \mathbb{N}^*} |y_n|\). Then \(e_n^*(x) = x_n\), and we have:
\[
\sum_{n=1}^\infty |e_n^*(x)| \, |y^*(e_n)| = \sum_{n=1}^\infty |x_n| \, |y_n|.
\]
Since \(|y_n| \leq \|y^*\|_{\ell_1^*}\) for all \(n\), it follows that:
\[
\sum_{n=1}^\infty |x_n| \, |y_n| \leq \|y^*\|_{\ell_1^*} \sum_{n=1}^\infty |x_n| = \|x\|_{\ell_1} \, \|y^*\|_{\ell_1^*}.
\]
Thus $\left( \left(e_{n},e_{n}^{* }\right) \right) _{n\in \mathbb{N}^{*}}$ is a Besselian frame for \(\ell_1\).\qed
\begin{theorem} \cite{kabbaj.karkri.zoubeir.2023}
\label{F.bess-equi.F*.bess}
A sequence $\mathcal{F}$ of $E$  is  Besselian
 if and only if 
 $\mathcal{F}^{*}$
 is a Besselian sequence of $E^{*}$. Moreover, in this case we have
$\mathcal{L}_{\mathcal{F}}=\mathcal{L}_{\mathcal{F}^{*}}$.
\end{theorem}
\begin{remark}
There is a Besselian sequence which is not a frame.
\end{remark}
\textbf{Proof.}
By theorem \eqref{F.bess-equi.F*.bess}, the dual sequence $\mathcal{F}^{*}$ of the Besselian frame from example \eqref{example.bess.frame} is a Besselian sequence in 
$\ell_{\infty}$. However, $\ell_{\infty}$  has no frames (as it is non-separable), and consequently  $\mathcal{F}^{*}$ cannot be a frame for $\ell_{\infty}$.\qed\\
\begin{lemma}
\label{Characterization.UCV}
 Let \(\sum_{m=1}^\infty z_m\) be a series in a Banach space \(E\). The following statements are equivalent:
 \begin{enumerate}
 \item
The series converges unconditionally.
\item
For every \(\varepsilon > 0\), there exists \(M \in \mathbb{N}^{*}\) such that for all finite sets \(A \subseteq \{M+1, M+2, \ldots\}\) and all scalars \(\lambda_m \in \mathbb{K}\) with \(|\lambda_m| \leq 1\):
   \[
   \left\| \sum_{m \in A} \lambda_m z_m \right\|_E < \varepsilon.
   \]
 \end{enumerate}
\end{lemma}

\textbf{Proof.} 
The implication $(2)\Rightarrow (1)$ is classical. Indeed, if condition (2) holds, then in particular (by taking $\lambda_{m}=1$ for all $m$) the series satisfies the Cauchy criterion and is therefore convergent. The unconditional convergence then follows from a standard characterization \cite[Theorem 16.1, p. 461]{I.Singer.I}.\\
 Now, assume that the series converges unconditionally.  Suppose, for contradiction, that condition (2) fails. Then there exists \(\varepsilon_0 > 0\) such that for every \(M\in \mathbb{N}^{*}\), there exist a finite set \(A_M \subseteq \{M+1, M+2, \ldots\}\) and scalars \((\lambda_m^M)_{m\in A_{M}}\) with \(|\lambda_m^M| \leq 1\) satisfying:
\[
\left\| \sum_{m \in A_M} \lambda_m^M z_m \right\|_E \geq \varepsilon_0.
\]

We  now construct a bounded sequence $(\lambda_m)_{m\in\mathbb{N}^{*}}$ such that the series \(\sum_{m=1}^{\infty}\lambda_m z_m\) diverges, contradicting the classical characterization of unconditional convergence.
\begin{itemize}
\item
Let \(M_1 = 1\). Choose a finite set \(A_{1} \subseteq \{M_1+1, M_1+2, \ldots\}\) and scalars \((\lambda_m^{1})_{m\in A_{1}}\) with \(|\lambda_m^{1}| \leq 1\) such that
$$\left\| \sum_{m \in A_{1}} \lambda_m^{1} z_m \right\|_E \geq \varepsilon_0.$$
\item
Set \(M_2 = \max A_{1} + 1\). Choose a finite set \(A_{2} \subseteq \{M_2+1, M_2+2, \ldots\}\) and scalars \((\lambda_m^{2})_{m\in A_{2}}\) with \(|\lambda_m^{2}| \leq 1\) such that
$$\left\| \sum_{m \in A_{2}} \lambda_m^{2} z_m \right\|_E \geq \varepsilon_0.$$
\item
Proceed inductively: having defined \(M_k\) and \(A_k\), set \(M_{k+1} = \max A_{k} + 1\), and choose a finite set \(A_{k+1} \subseteq \{M_{k+1}+1, M_{k+1}+2,\ldots\}\) and scalars \((\lambda_m^{k+1})_{m\in A_{k+1}}\) with \(|\lambda_m^{k+1}| \leq 1\) such that
$$\left\| \sum_{m \in A_{k+1}} \lambda_m^{k+1} z_m \right\|_E \geq \varepsilon_0.$$
 \end{itemize}
Define the sequence \((\lambda_m)_{m\in\mathbb{N}^{*}}\) by:
\[
\lambda_m = 
\begin{cases}
\lambda_m^{k} & \text{if } m \in A_{k} \text{ for some } k, \\
0 & \text{otherwise}.
\end{cases}
\]
By construction, \(|\lambda_m| \leq 1\) for all $m$, so \((\lambda_m)_{m\in\mathbb{N}^{*}}\) is bounded. Now consider the partial sums over the finite sets \(S_k = \bigcup_{j=1}^k A_{j}\). The sets $(A_k)_{k\in\mathbb{N}^{*}}$ are finite, disjoint, and their union is increasing. For each $k\geq 2$, we have:
\[
\left\| \sum_{m \in S_k} \lambda_m z_m - \sum_{m \in S_{k-1}} \lambda_m z_m \right\|_E = \left\| \sum_{m \in A_{k}} \lambda_m^{k} z_m \right\|_E \geq \varepsilon_0.
\]
This shows that the sequence of partial sums $\left(\sum_{m\in S_{k}} \lambda_m z_m\right)_{k\in\mathbb{N}^{*}}$ is not Cauchy. Therefore, the series $\sum \lambda_m z_m$ diverges, contradicting the assumption of unconditional convergence.\qed\\

The following theorem is inspired by  \cite[p. 32-33]{A.Abramovich.2002}.
\begin{theorem}
\label{unc.impl.bess}
Suppose that $\mathcal{F}$ is an unconditional frame for $E$. Then $\mathcal{F}$ is a Besselian frame for $E$.
\end{theorem}
\textbf{Proof}.
For each  $M\in \mathbb{N}^{*}$  and 
 $(s_{m})_{m\in\mathbb{N}^{*}}\in\mathcal{S}$, the linear operator
  $$T_{M,s}:E\rightarrow E,\; x\mapsto\sum_{m=1}^{M}s_{m}b_{m}^{*}(x)a_{m}$$ 
is  bounded. By unconditionality of $\mathcal{F}$ and the Banach-Steinhaus theorem, the linear operator
  $$T_{s}:E\rightarrow E,\; x\mapsto\sum_{m=1}^{\infty} 
  s_{m}b_{m}^{*}(x)a_{m}$$ is bounded.
We now show that for each $x\in E$, the mapping 
 $s\mapsto T_{s}(x) $ is continuous on the compact metric space
 $\mathcal{S}$.\\
  Fix $x\in E$ and $\varepsilon > 0$. By lemma \eqref{Characterization.UCV}, there exists $M\in\mathbb{N}^{*}$ such that for all finite sets 
 $A \subseteq\{M+1, M+2, \dots\}$ and all scalars $\lambda_{m}\in\mathbb{K}$ with 
 $|\lambda_m| \leq 1$:  
  \[
  \left\|\sum_{m\in A}\lambda_m b_m^*(x) a_m \right\|_E < \frac{\varepsilon}{4}.
  \]  
Now, for \(s, t \in \mathcal{S}\), we estimate the difference:
  \[
  \|T_s(x) - T_t(x)\|_E \leq \underbrace{\left\| \sum_{m=1}^M(s_m - t_m) b_m^*(x) a_m \right\|_E}_{\text{(I)}} + \underbrace{\left\| \sum_{m=M+1}^\infty (s_m - t_m) b_m^*(x) a_m \right\|_E}_{\text{(II)}}.
  \]  
Estimate (II): Since \(|s_m - t_m| \leq 2\), define \(\lambda_m = \frac{s_m - t_m}{2}\) (\(|\lambda_m| \leq 1\)). Then,
    \[
    \text{(II)} \leq 2  \sup_{A \subseteq \{M+1,\dots\}} \left\| \sum_{m \in A} \lambda_m b_m^*(x) a_m \right\|_E \leq 2  \frac{\varepsilon}{4} = \frac{\varepsilon}{2}.
    \]  
 Estimate (I): Let \(C = \max_{1 \leq m \leq M} \|b_m^*(x) a_m\|_E\). Choose 
   $$\delta = \frac{\varepsilon}{2^{M+1}MC}.$$
If \(d(s,t) < \delta\), then for each \(m \leq M\),
$$|s_m - t_m|<2^m d(s,t)< 2^m \delta \leq 2^M\delta .$$
Hence,
  \[
    \text{(I)} \leq \sum_{m=1}^M |s_m - t_m|  C <M2^M\delta C = \frac{\varepsilon}{2}.
    \]  
  Combining (I) and (II), we conclude that:
$$\|T_s(x) - T_t(x)\|_E < \varepsilon$$
 whenever \(d(s,t) < \delta\). Thus, $s\mapsto T_{s}(x) $ is continuous on 
 $\mathcal{S}$.\\
  Since $\mathcal{S}$ is compact, the set
 $\left\lbrace \Vert T_{s}(x)\Vert_{E}: s\in \mathcal{S}\right\rbrace$ is bounded. 
By the uniform boundedness principle,  there is a constant $\mathcal{D}>0$ such that
$$\left \Vert \underset{m=1}{\overset{\infty}{\sum}}
s_{m}b_{m}^{*}(x)a_{m}\right\Vert_{E}
\leq
\mathcal{D}\Vert x\Vert_{E},
\;\;\;\left (x\in E,\;(s_{m})_{m\in\mathbb{N}^{*}}\in \mathcal{S}\right ).$$
Now fix $x\in E$ and $y^{*}\in E^{*}$. Choose a sequence of signs
$s=\left( s_{m}\right)_{m\in \mathbb{N}^{*}}$ such that 
$$\vert b_{m}^{*}(x)y^{*}(a_{m})\vert=s_{m}b_{m}^{*}(x)y^{*}(a_{m}),\;\;\;(m\in \mathbb{N}^{*}).$$
(with $s_{m}=1$ if $b_{m}^{*}(x)y^{*}(a_{m})=0 $) Then:
\begin{align*}
\underset{m=1}{\overset{\infty}{\sum}}
\left\vert b_{m}^{*}(x)y^{*}(a_{m})\right\vert 
& =\overset{\infty}{\underset{m=1}{\sum }}
s_{m}b_{m}^{*}(x)y^{*}(a_{m})
=\left\vert y^{*}\left(\overset{\infty }{\underset{m=1}{\sum}}
s_{m}b_{m}^{*}(x)a_{m}\right)\right\vert\\
&=\left\vert y^{*}\left(T_{s}(x)\right)\right\vert 
\leq
\Vert y^{*}\Vert_{E^{*}}\left\Vert T_{s}(x)\right\Vert_{E}\\
&\leq
\Vert y^{*}\Vert_{E^{*}}\left\Vert T_{s}\right\Vert_{L(E)}\cdot\left\Vert x\right\Vert_{E}\\
&\leq \mathcal{D}\left\Vert x\right\Vert_{E}\Vert y^{*}\Vert_{E^{*}}.
\end{align*}
Thus $\mathcal{F}$ is a Besselian Schauder frame for $E$.\qed\\   

\begin{example}
For every $p>1$, the classical Banach spaces
$\ell_{p}$ and $L_{p}([0,1])$ have Besselian frames.\qed
\end{example}
\textbf{Proof.} The canonical basis of $\ell_{p}$ (resp. of $L_{p}([0,1]))$  is unconditional. By theorem \eqref{unc.impl.bess}, it follows that this basis is a Besselian frame for $\ell_{p}$ (resp. for $L_{p}([0,1]))$.\qed
\begin{corollary}
\label{unconditional.basis.equi.bess}
Suppose that $\mathcal{F}$ is a basis for $E$. Then $\mathcal{F}$
is unconditional if and only if it is Besselian.
\end{corollary}
\textbf{Proof.}
The equivalence follows immediately from \cite[theorem 16.1, p. 461]{I.Singer.I} and theorem \eqref{unc.impl.bess}.\qed
\begin{corollary}
There is a frame which is not Besselian.
\end{corollary}
\textbf{Proof.}
It is well-known \cite[example 3.1.2, p. 51]{Albiac.Kalton.2006} that the sequence
$$\left( \left(\sum_{k=1}^{n}e_{k},e_{n}^{*}-e_{n+1}^{*}\right) \right) _{n\in \mathbb{N}^{*}}$$
forms a conditional (i.e., not unconditional) basis (and thus a frame) for 
$c_{0}$. By corollary \eqref{unconditional.basis.equi.bess}, this frame is not Besselian.\qed
\begin{corollary}
\label{wsc.then.unc.equi.bess}
Suppose that $E$ is WSC and $\mathcal{F}$ is a frame for $E$. Then $\mathcal{F}$
is unconditional if and only if it is Besselian.
\end{corollary}
\textbf{Proof.}
If $\mathcal{F}$ is unconditional.
Theorem \eqref{unc.impl.bess} implies that $\mathcal{F}$ is Besselian. Conversely, if $\mathcal{F}$ is Besselian, then for every $x\in E$, the series $\sum _{n}b_{n}^{*}(x)a_{n}$ is  WUC in the WSC Banach space $E$.
By \cite[corollary 2.4.15, p. 43]{Albiac.Kalton.2006},  this series $\sum _{n}b_{n}^{*}(x)a_{n}$ converges unconditionally in
$E$. Thus, $\mathcal{F}$ is unconditional. \qed
\begin{corollary}
\label{L_1.has.no.Bess}
The space $L_{1}([0,1])$ has no Besselian frame.
\end{corollary}
\textbf{Proof.} Assume, by contradiction, that $L_{1}([0,1])$ has a Besselian frame $\mathcal{F}$. 
By corollary \eqref{wsc.then.unc.equi.bess}, $\mathcal{F}$
would then be an unconditional frame for $L_{1}([0,1])$. According to \cite[proposition 2.4, p. 69]{Casazza.D.S.Z.2008},  this would imply that $L_{1}([0,1])$ is embedded in a Banach space with unconditional basis, which  contradicts the fact that $L_{1}([0,1])$ cannot be embedded in a Banach space with unconditional basis \cite[theorem 6.3.3, p. 144]{Albiac.Kalton.2006}.\qed
\begin{proposition}
\label{injective.projective.tensor.no.bess}
For each  $p>q>1$, $\ell_{p^{*}}\otimes_{\varepsilon}\ell_{q}$ 
and $\ell_{p}\otimes_{\pi}\ell_{q^{*}}$
have no Besselian  frames (then have no unconditional frame).
\end{proposition}
\textbf{Proof.}
Suppose that $\mathcal{F}$ is a Besselian frame for  $\ell_{p^{*}}\otimes_{\varepsilon}\ell_{q}$ 
(resp.  $\ell_{p}\otimes_{\pi}\ell_{q^{*}}$). Since
$\ell_{p^{*}}\otimes_{\varepsilon}\ell_{q}$ 
(resp. $\ell_{p}\otimes_{\pi}\ell_{q^{*}}$) 
is a reflexive space \cite[corollary 4.24, p. 87]{Raymond.A.Ryan} then it is  WSC. 
It follows by corollary \eqref{wsc.then.unc.equi.bess} that $\mathcal{F}$ is unconditional.
Therefore $\ell_{p^{*}}\otimes_{\varepsilon}\ell_{q}$ 
(resp. $\ell_{p}\otimes_{\pi}\ell_{q^{*}}$) is complemented in a Banach space with an unconditional basis \cite[proposition 2.4, p. 69]{Casazza.D.S.Z.2008}, which  contradicts the fact that 
$\ell_{p^{*}}\otimes_{\varepsilon}\ell_{q}$ 
(resp.  $\ell_{p}\otimes_{\pi}\ell_{q^{*}}$) is not isomorphic to a subspace 
of a Banach space with an unconditional basis \cite[corollary 3.1, p. 58]{Kwapien.Pelski.1970}.\qed
\begin{lemma}\cite[corollary 4.13, p. 80]{Raymond.A.Ryan}
\label{compact.injective.tensor}
If either $E^{*}$ or $F$ has the approximation property, then
$K(E,F)=E^{*}\otimes_{\varepsilon}F$.
\end{lemma}
\begin{corollary}
For each  $p>q>1$, $K(\ell_{p},\ell_{q})$ 
has no Besselian  frame (then has no unconditional frame).
\end{corollary}
\textbf{Proof.}
It follows directly from lemma \eqref{compact.injective.tensor}
and proposition \eqref{injective.projective.tensor.no.bess}.\qed
\begin{lemma}
\label{WUV.not.UCV}
There exists a series \(\sum x_n\) in \(c_0\) such that:  
$$\sum_{n=1}^\infty \vert y^{*}(x_n)\vert \leq  2\|y^{*}\|_{\ell_{1}} ,\;\;(y^{*} \in c_0^{*})$$ 
and \(\sum_{n=1}^\infty x_n\) converge conditionally (i.e., it converges but not unconditionally) in \(c_0\) to $0$.
\end{lemma}   
\textbf{Proof.}
We construct a series \(\sum x_n\) in \(c_0\) as follows. For each \(k \ge 1\), define a block of \(2k\) terms, all affecting only the \(k\)-th coordinate.
 For \(j = 1, \dots, k\):  
  $$x_{N_k + 2j - 1} = \dfrac{1}{k} e_k,$$  
 $$x_{N_k + 2j} = -\dfrac{1}{k} e_k,$$
where \(N_k = 2 \sum_{i=1}^{k-1} i = k(k-1)\) (with \(N_1 = 0\)), and \((e_k)_{k\in\mathbb{N}^{*}}\) is the standard unit vector basis of \(c_0\). The terms are ordered by completing each block consecutively.
\begin{enumerate}
\item
Convergence in \(c_0\):\\
 The partial sum after each complete block is the zero vector. Within block \(k\), the sup-norm of any partial sum is at most \(1/k\). Given \(\varepsilon > 0\), choose \(K\) such that \(1/K < \varepsilon\). For all \(n > N_K\), the partial sum \(S_n\) involves only blocks with index \(\ge K+1\) or an incomplete blocks of index \(\le K\). In either case \(\|S_n\|_{c_{0}} \le 1/K < \varepsilon\). Hence \(\sum x_n\) converges to \(0\) in \(c_0\).
\item
Weakly unconditionally Cauchy (WUC):\\ 
  The dual space \(c_0^{*}\) is isometrically isomorphic to \(\ell_1\).  
   Every \(y^{*} \in c_0^{*}\) corresponds to a sequence \((y^{*}_k) \in \ell_1\) with \(\|y^{*}\|_{\ell_1} = \sum_{k=1}^{\infty} |y^{*}_k|\).  
   For each \(k\), the block contributes \(2k\) terms, each of absolute value \(|y^{*}(x_n)| = |y^{*}_k|/k\).  
   Therefore, 
   \[
   \sum_{n=1}^\infty |y^{*}(x_n)| = \sum_{k=1}^\infty 2k \cdot \frac{|y^{*}_k|}{k} = 2 \sum_{k=1}^\infty |y^{*}_k| = 2\|y^{*}\|_{\ell_{1}} < \infty.
   \]
   Thus the series is WUC.
\item
Not unconditionally convergent:\\ 
A series \(\sum x_n\) in a Banach space converges unconditionally if and only if \(\sum \varepsilon_n x_n\) converges for every choice of signs \(\varepsilon_n \in \{\pm 1\}\).
 Choose \(\varepsilon_n\) as follows:  
   \[
   \varepsilon_n = 
   \begin{cases}
   +1 & \text{if } x_n = \frac{1}{k} e_k, \\
   -1 & \text{if } x_n = -\frac{1}{k} e_k.
   \end{cases}
   \]
 Hence,
   \[
   \sum \varepsilon_n x_n = \sum_{k=1}^{\infty} (2k \text{ terms of } \frac{1}{k} e_k) = \sum_{k=1}^{\infty} 2 e_k.
   \]

   The partial sums are \(S_K = \sum_{k=1}^{K} 2 e_k\), and
   \[
   \|S_{K+1} - S_K\|_{c_{0}} = \|2 e_{K+1}\|_{c_{0}} = 2,
   \]
   which does not tend to \(0\). Thus the series \(\sum \varepsilon_n x_n\) does not converge in \(c_0\).  
   By the sign criterion, the original series \(\sum x_n\) is not unconditionally convergent.
\end{enumerate}
Thus, \(\sum x_n\) is conditionally convergent (not unconditionally) and weakly unconditionally Cauchy in \(c_0\).
 \begin{theorem}
 \label{Besselian.not.unconditional}
There exists a Besselian frame which is not unconditional.
\end{theorem}

\textbf{Proof}\\
We construct a frame $(a_n, b^{*}_{n})_{n\in\mathbb{N}^{*}}$ for \(c_0\) using the series from  lemma \eqref{WUV.not.UCV}. Let $\sum x_n$  be the series in \(c_0\) constructed in lemma \eqref{WUV.not.UCV}, which converges conditionally to \(0\) and is weakly unconditionally Cauchy (WUC). Let \((e_k)_{k\in\mathbb{N}^{*}}\) denote the standard unit vector basis of \(c_0\) and \((e^{*}_k)_{k\in\mathbb{N}^{*}}\) the corresponding coordinate functionals. Fix a nonzero functional \(\varphi \in c_0^*\); for instance, take \(\varphi = e_1^*\).

Define the frame by interleaving terms as follows. For each 
\(k = 1, 2, 3, \dots\):
  \[
  a_{2k-1} = e_k, \quad b_{2k-1}^* = e_k^*.
  \]
  \[
  a_{2k} = x_k, \quad b_{2k}^* = \varphi.
  \]

We verify the required properties:
\begin{enumerate}
\item
Frame property:\\
 For any \(x \in c_0\),
   \[
   \sum_{n=1}^\infty b_n^*(x) a_n = \sum_{k=1}^\infty e_k^*(x) e_k + \varphi(x) \sum_{k=1}^\infty x_k = x + \varphi(x) \cdot 0 = x,
   \]
   since \(\sum_{k=1}^\infty x_k = 0\) by lemma \eqref{WUV.not.UCV}. Thus $(a_n, b_n^*)_{n\in\mathbb{N}^{*}}$ is a frame for $c_{0}$.
\item
Besselian Condition:\\
 For any \(x \in c_0\) and \(y^* \in c_0^* \cong \ell^1\), using  example \eqref{example.bess.frame.c_0} and lemma\eqref{WUV.not.UCV} we obtain
  \begin{align*}  
   \sum_{n=1}^\infty |b_n^*(x) y^*(a_n)| 
  &= \sum_{k=1}^\infty |e_k^*(x) y^*(e_k)| + \sum_{k=1}^\infty |\varphi(x) y^*(x_k)|\\
& \leq
\|x\|_{c_{0}} \|y^{*}\|_{\ell_{1}}+\|x\|_{c_{0}}\|\varphi\|_{\ell_{1}}\sum_{k=1}^\infty |y^*(x_k)|\\
 &= 
 3\|x\|_{c_{0}} \|y^{*}\|_{\ell_{1}}
\end{align*}
 which establishes that $(a_n, b_n^*)_{n\in\mathbb{N}^{*}}$ is a Besselian frame with constant \(\mathcal{L}_{\mathcal{F}} \leq 3\).
\item
Conditional convergence for some \(x\):\\
 Take \(x = e_1\). Then $\varphi(e_1)=1$, and the frame expansion becomes
   \[
   \sum_{n=1}^\infty b_n^*(e_1) a_n = e_1 + \sum_{k=1}^\infty x_k.
   \]
  Since \(\sum_{k=1}^{\infty} x_k = 0\), this series converges to \(e_1\) in the given ordering.  
   However, because \(\sum_{k=1}^{\infty} x_k\) converges conditionally to \(0\), there exists a rearrangement \(\sum_{k=1}^{\infty} x_{\sigma(k)}\) that does not converge to \(0\) (or converges to a different limit). Consequently, rearranging the even-indexed terms of the frame expansion accordingly yields a series that does not converge to \(e_1\). Therefore, the series \(\sum_{n=1}^{\infty} b_n^*(e_1) a_n\) does not converge unconditionally to \(e_1\), i.e., it converges conditionally.
\end{enumerate}
Thus, we have constructed a Besselian frame for \(c_0\) and an element \(e_1 \in c_0\) such that the frame expansion converges conditionally to \(e_1\). 
\begin{theorem}
\label{zabr}
$\mathcal{F}$ is a Besselian sequence of $E$  if and only if
for every $x\in E$ and $y^{*}\in E^{*}$, the following condition holds:
\begin{equation}
\underset{n=1}{\overset{\infty }{\sum }}
\left \vert b_{n}^{*}\left(x\right) \right \vert \left \vert y^{*}\left(a_{n}\right)
\right \vert <\infty. \label{bess.characterization}
\end{equation}
\end{theorem}
\textbf{Proof.}
Suppose that the inequality \eqref{bess.characterization} holds
for each $x\in E$ and $y^{*}\in E^{*}$.
Let $y^{*}\in E^{*}$ and   
 $(s_{k})_{k\in\mathbb{N}^{*}}\in \mathcal{S}$, 
an argument analogous to theorem \eqref{unc.impl.bess} yields a constant 
 $\mathcal{D}_{y^{*}}>0$ such that:
$$\left \vert \underset{n=1}{\overset{\infty }{\sum }}
s_{n}b_{n}^{*}(x)y^{*}(a_{n}) \right \vert \leq \mathcal{D}_{y^{*}}\Vert x\Vert_{E},\;\;\;\left (x\in E,\;(s_{n})_{n\in\mathbb{N}^{*}}\in \mathcal{S}\right ).$$
Consequently, we obtain:
\begin{align}
\label{zabreiko.1}
\underset{n=1}{\overset{\infty }{\sum }}
\left \vert b_{n}^{*}(x)y^{*}(a_{n})\right\vert 
\leq \mathcal{D}_{y^{*}}\left \Vert  x \right \Vert_{E}, 
\;\;\; (x\in E, y^{*}\in E^{*}).
\end{align}
Similarly, for each $x\in E$ there is a constant
 $\mathcal{C}_{x}>0$ such that
 $$\left \vert \underset{n=1}{\overset{\infty }{\sum }}
s_{n}b_{n}^{*}(x)y^{*}(a_{n}) \right \vert \leq \mathcal{C}_{x}\Vert y^{*}\Vert_{E^{*}},\;\;\;\left (y^{*}\in E^{*},\;(s_{n})_{n\in\mathbb{N}^{*}}\in \mathcal{S}\right ),$$
therefore we obtain:
\begin{align}
\label{zabreiko.2}
\underset{n=1}{\overset{\infty }{\sum }}
\left \vert b_{n}^{*}(x)y^{*}(a_{n})\right\vert 
\leq \mathcal{C}_{x}\left \Vert y^{*} \right \Vert_{E^{*}}, 
\;\;\; (x\in E, y^{*}\in E^{*}).
\end{align}
Finally, combining inequalities \eqref{zabreiko.2} and 
\eqref{zabreiko.1} with the method of \cite[theorem 3.6, p. 232]{kabbaj.karkri.zoubeir.2023} we conclude that $\mathcal{F}$ is Besselian. The reverse implication is immediate from the definition.\qed\\

The existence of Banach spaces with unconditional frames that lack unconditional bases remains an open problem. In other words, is the inclusion
\[
\left\lbrace\text{Banach spaces with}\; USB\right\rbrace \subsetneq \{\text{Banach spaces with}\;USF\}
\]
strict?  (see \cite[remark 5.8, p. 170]{Casazza.D.H.L.1999} ).\\

\textbf{Problem.} Does there exist a Banach space \(E\) that admits a Besselian frame but no unconditional frame? In other words, is the inclusion
\[
\left\lbrace\text{Banach spaces with}\; USF\right\rbrace \subsetneq \{\text{Banach spaces with}\;BSF\}
\]
strict? A positive answer,  together with corollary \eqref{L_1.has.no.Bess}  would show that the class of spaces with Besselian frames lies strictly between the class of spaces with unconditional frames and the class of spaces with frames.  These two problems form a natural hierarchy of questions about the strictness of the inclusions between classes of Banach spaces defined by increasingly general frame properties.
\section{Structures of E, E* and E**}
In this section, we generalize some  well-known results of Karlin \cite{S.Karlin.1948}  and  James \cite{James.1950}
to Banach spaces with Besselian  frames (then in particular
  with unconditional  frames).
 \begin{proposition}
If $E^{*}$ has a frame, then $E$ has a frame.
\end{proposition}
\textbf{Proof.}
Suppose that $E^{*}$ has a  frame. Then  $E^{*}$ has the BAP, and by \cite[theorem 4.9, p. 294]{Casazza.AP.2001} $E$ is embedded complementably into a space with a shrinking basis, consequently
 $E$ has a frame.\qed
\begin{proposition} 
There is a separable Banach space with a frame whose dual is
separable but fails to have a  frame.
\end{proposition}
\textbf{Proof.}
By \cite[theorem 1.e.7, p. 34]{Lindenstraus.Tzafriri.1977},
there is a Banach space $E$ with a basis (and consequently a frame) such that $E^{*}$ is separable but fails to have the AP. Then, $E^{*}$ fails to have the BAP. Therefore, $E^{*}$ fails to have a frame.\qed\\

In this section, we will make extensive use of the following lemma
(which is  inspired by \cite{S.Karlin.1948}).
\begin{lemma}
\label{interchange.E}
Suppose that $\mathcal{F}$ is a Besselian frame for $E$, $\left(x_{r}\right)_{r\in\mathbb{N}^{*}}$ is a weakly Cauchy sequence 
 in $E$ and $y^{*}\in E^{*}$. Then we have
\begin{align}
\underset{r\rightarrow \infty}{\lim }
y^{*}(x_{r})
=
 \underset{r\rightarrow \infty}{\lim }
 \underset{n=1}{\overset{\infty }{\sum }}
  b_{n}^{*}(x_{r})y^{*}\left(a_{n}\right)
=
\underset{n=1}{\overset{\infty}{\sum }}
\underset{r\rightarrow \infty}{\lim }
 b_{n}^{*}(x_{r})y^{*}\left(a_{n}\right).
\end{align}
\end{lemma}
\textbf{Proof.} 
Since $\mathcal{F}$ is Besselian, the first equality follows immediately from the frame representation 
$$x_{r}
=
 \underset{n=1}{\overset{\infty }{\sum }}
  b_{n}^{*}(x_{r})a_{n},\;\;\left(r\in\mathbb{N}^{*}\right).$$
Consider the operator:
 $$U:E\rightarrow \ell_{1};\;\;
 x\mapsto \left(  b_{n}^{*}(x)y^{*}\left(a_{n}\right)\right)_{n\in\mathbb{N}^{*}}.$$ 
The Besselian property ensures that $U$ is bounded. Thus, 
 $\left(  U(x_{r})\right)_{r\in\mathbb{N}^{*}}$ is a weakly Cauchy sequence in $\ell_{1}$. By the Schur property of $\ell_{1}$, this implies that $\left(  U(x_{r})\right)_{r\in\mathbb{N}^{*}}$ is a Cauchy sequence.
 Let $\varepsilon>0$. There exists $r_{0}\in\mathbb{N}^{*}$ such that:
 $$\underset{n=m}{\overset{\infty }{\sum }}
 \left \vert b_{n}^{*}(x_{r})y^{*}\left(a_{n}\right)
 -b_{n}^{*}(x_{r_{0}})y^{*}\left(a_{n}\right)   \right \vert
\leq \varepsilon, \;\;(r\geq r_{0},\;\;m\in\mathbb{N}^{*}).$$
Since the series $\sum_{n}b_{n}^{*}(x_{r_{0}})y^{*}\left(a_{n}\right)$
is convergent, there exists $N_{0}\in\mathbb{N}^{*}$ such that
 for each $m\geq N_{0}$ we have
$\left \vert \sum_{m}^{\infty} b_{n}^{*}(x_{r_{0}})y_{k}^{*}\left(a_{n}\right)
    \right \vert\leq \varepsilon$. Hence
 \begin{align*}
 \left \vert \underset{n=1}{\overset{\infty }{\sum }}
  b_{n}^{*}(x_{r})y^{*}\left(a_{n}\right)
 -\underset{n=1}{\overset{m}{\sum }}
 b_{n}^{*}(x_{r})y^{*}\left(a_{n}\right)   \right \vert
&\leq 2\varepsilon, \;\;\;(r\geq r_{0},\;m\geq N_{0}).
 \end{align*}
 Then
$$
 \left \vert \underset{r\rightarrow \infty}{\lim }
 \underset{n=1}{\overset{\infty }{\sum }}
  b_{n}^{*}(x_{r})y^{*}\left(a_{n}\right)
 -\underset{n=1}{\overset{m}{\sum }}
 \underset{r\rightarrow \infty}{\lim }
  b_{n}^{*}(x_{r})y^{*}\left(a_{n}\right)   \right \vert
\leq 2\varepsilon, \;\;\;(m\geq N_{0}),
$$
and consequently we have
\begin{align}
\label{karlin.1}
 \underset{r\rightarrow \infty}{\lim }
 \underset{n=1}{\overset{\infty }{\sum }}
  b_{n}^{*}(x_{r})y^{*}\left(a_{n}\right)
=
\underset{n=1}{\overset{\infty}{\sum }}
\underset{r\rightarrow \infty}{\lim }
 b_{n}^{*}(x_{r})y^{*}\left(a_{n}\right).
\end{align}
Which finishes the proof.\qed
 \begin{lemma}
 \label{interchange.E*}
Suppose that   $\mathcal{F}$ is a Besselian frame for $E$, $\left(y_{k}^{*}\right)_{k\in\mathbb{N}^{*}}$ is a weakly Cauchy sequence 
in $E^{*}$ and  $z^{**}\in E^{**}$.  If the limit 
 of $\sum_{n=1}^{\infty}
  z^{**}(b_{n}^{*})y_{k}^{*}\left(a_{n}\right) $ exists as $k$ tends to infinity, then:
\begin{align}
\label{karlin.1}
\underset{k\rightarrow \infty}{\lim }
 \underset{n=1}{\overset{\infty }{\sum }}
  z^{**}(b_{n}^{*})y_{k}^{*}\left(a_{n}\right)
=
\underset{n=1}{\overset{\infty}{\sum }}
\underset{k\rightarrow \infty}{\lim }
 z^{**}(b_{n}^{*})y_{k}^{*}\left(a_{n}\right).
\end{align}
\end{lemma}
\textbf{Proof.}
Since $\mathcal{F}$ is Besselian, its dual sequence 
$\mathcal{F}^{*}$ is also Besselian. Consider the operator:
$$V:E^{*}\rightarrow \ell_{1};\;\;
y^{*}\mapsto \left( z^{**}(b_{n}^{*})y^{*}\left(a_{n}\right)\right)_{n\in\mathbb{N}^{*}}.$$ 
The Besselian property ensures that $V$ is bounded. Thus, 
 $\left( V(y_{k}^{*})\right)_{k\in\mathbb{N}^{*}}$ is a weakly Cauchy sequence, and hence 
 is a Cauchy sequence.
The conclusion now follows by an argument analogous to that in lemma \eqref{interchange.E}.\qed
\begin{lemma}
 \label{interchange.E**}
Let $\mathcal{G}:=\left( \left( b_{n}^{*},\phi_{n}^{**}\right) \right)
_{n\in \mathbb{N}^{*}}$ be a Besselian frame for $E^{*}$.
Let $\left(x_{r}\right)_{r\in\mathbb{N}^{*}}$ be a weakly Cauchy sequence 
in $E$ and $z^{***}\in E^{***}$. If the limit 
 of:
$$\sum_{n=1}^{\infty}
  z^{***}(\phi_{n}^{**})b_{n}^{*}(x_{r}),$$
exists as $r$ tends to infinity. Then
\begin{align}
\underset{r\rightarrow \infty}{\lim }
 \underset{n=1}{\overset{\infty }{\sum }}
  z^{***}(\phi_{n}^{**})b_{n}^{*}\left(x_{r}\right)
=
\underset{n=1}{\overset{\infty}{\sum }}
\underset{r\rightarrow \infty}{\lim }
 z^{***}(\phi_{n}^{**})b_{n}^{*}\left(x_{r}\right).
\end{align}
\end{lemma}
\textbf{Proof.}
Since $\mathcal{G}$ is Besselian, its dual sequence $\mathcal{G}^{*}$ is also Besselian. Consider the operator:
$$W:E\rightarrow \ell_{1};\;\;
 x\mapsto \left(z^{***}(\phi_{n}^{**})(b_{n}^{*}(x)\right)_{n\in\mathbb{N}^{*}}.$$
The Besselian property ensures that $W$ is bounded. Thus,
 $\left(W(x_{r})\right)_{r\in\mathbb{N}^{*}}$ is a weakly Cauchy sequence in $\ell_{1}$ and hence is a  Cauchy sequence. The conclusion now follows by an argument analogous to that in lemma \eqref{interchange.E}.
\qed\\ 

The following theorem extends  the Karlin theorem \cite[p. 981]{S.Karlin.1948} to Banach spaces with Besselian frames.
\begin{theorem}
\label{bess.dual.WSC}
Suppose that  $\mathcal{F}$ is a Besselian frame for  $E$ and $E^{*}$ is separable. Then $E^{*}$ is  weakly sequentially complete.
\end{theorem}
\textbf{Proof.}
Let $\left( y_{k}^{*}\right)_{k\in\mathbb{N}^{*}}$ be a weakly Cauchy
 sequence in $E^{*}$.
Then, for each $\phi\in E^{**}$  the limit of $\phi(y_{k}^{*})$ exists as $k$ tends to infinity.
In particular, there exists an element $y_{0}^{*}\in E^{*}$ such that $y_{k}^{*}(x)\rightarrow y_{0}^{*}(x)$
for each $x\in E$. To prove the theorem, it suffices to show that
$\phi(y_{k}^{*})\rightarrow \phi(y_{0}^{*})$ for every $\phi\in E^{**}$. 
Indeed, let $\phi\in E^{**}$. Since $E^{*}$ is separable, Gelfand theorem (see \cite{Gantmakhter.Smulian.1937}) guarantees the existence of a sequence 
$\left(x_{r}\right)_{r\in\mathbb{N}^{*}}$ in $E$ such that:
\begin{align}
\label{Gelfand}
\phi(y^{*})=\underset{r\rightarrow \infty}{\lim }y^{*}(x_{r}),\;\;\;(y^{*}\in E^{*}).
\end{align}
Consequently, by lemmas \eqref{interchange.E} and \eqref{interchange.E*}, we obtain:
\begin{align*}   
\underset{k\rightarrow \infty}{\lim }\phi(y_{k}^{*})
&=\underset{k\rightarrow \infty}{\lim }\underset{r\rightarrow \infty}{\lim }y_{k}^{*}(x_{r})
=\underset{k\rightarrow \infty}{\lim }
\underset{r\rightarrow \infty}{\lim }\overset{\infty }{\underset{n=1}{\sum }}
 b_{n}^{*}(x_{r})y_{k}^{*}(a_{n})\\
 & =\underset{k\rightarrow \infty}{\lim }
 \overset{\infty }{\underset{n=1}{\sum }}
\phi( b_{n}^{*})y_{k}^{*}(a_{n})
=\overset{\infty }{\underset{n=1}{\sum }}
\phi( b_{n}^{*})y_{0}^{*}(a_{n})\\
 &=\underset{r\rightarrow \infty}{\lim }
 \overset{\infty }{\underset{n=1}{\sum }}
 b_{n}^{*}(x_{r})y_{0}^{*}(a_{n})
 =\underset{r\rightarrow \infty}{\lim }\underset{k\rightarrow \infty}{\lim }\overset{\infty }{\underset{n=1}{\sum }}
 b_{n}^{*}(x_{r})y_{k}^{*}(a_{n})\\
 &=\underset{r\rightarrow \infty}{\lim }\underset{k\rightarrow \infty}{\lim }y_{k}^{*}(x_{r})=\phi(y_{0}^{*}).
 \end{align*}
Therefore, $E$ is   weakly sequentially complete.\qed
 \begin{remark}
The argument \eqref{Gelfand} in the proof of Theorem \eqref{bess.dual.WSC} works assuming only that $E$ is separable and contains no isomorphic copy of $\ell_{1}$, since by Odell-Rosenthal theorem, in this case $E$ is   weak$^{*}-$sequentially dense in $E^{**}$.
 \end{remark}
\begin{lemma}\cite[proposition 3.3, p. 230]{kabbaj.karkri.zoubeir.2023}
\label{E*.wsc.bess}
Suppose that $\mathcal{F}$ is a Besselian frame for  $E$, and $E^{*}$ is 
WSC. Then $\mathcal{F^{*}}$ is a Besselian (then unconditional) frame for $E^{*}$.
\end{lemma}
\begin{corollary}
 The Banach space $C[0,1]$  has no Besselian  frame.
 \end{corollary}
\textbf{Proof.} Assume, for contradiction, that $C[0,1]$ has a Besselian frame. Since $C[0,1]^{*}=V$ (space of functions of bounded variation), lemma \eqref{E*.wsc.bess} implies that $V$ must also have a Besselian frame. However, this contradicts the non-separability of $V$.\qed\\

The following corollary extends  the Karlin corollary \cite[p. 983]{S.Karlin.1948} to Banach spaces with Besselian frames.
\begin{corollary}
\label{bess.dual.bess}
If $E$ has a Besselian  frame and its dual $E^{*}$ is separable, then
$E^{*}$ has a Besselian frame.
\end{corollary}
\textbf{Proof.} This follows immediately from theorem \eqref{bess.dual.WSC} and
 lemma \eqref{E*.wsc.bess}.\qed
 \begin{corollary}
The space  $K(\ell_{2},\ell_{2})$
has no Besselian frame (and consequently has no unconditional frame).
\end{corollary}
\textbf{Proof.}
Assume, for contradiction, that $K(\ell_{2},\ell_{2})$ has a Besselian frame $\mathcal{F}$. Since $\left (K(\ell_{2},\ell_{2})\right )^{*}=\ell_{2}\otimes_{\pi}\ell_{2}$ 
\cite[p. 91]{Raymond.A.Ryan}, theorem \eqref{bess.dual.WSC} implies that 
$\ell_{2}\otimes_{\pi}\ell_{2}$ is WSC. 
By corollary \eqref{wsc.then.unc.equi.bess}, it follows that $\mathcal{F}$ must be unconditional.
Consequently, $K(\ell_{2},\ell_{2})$ is complemented in a Banach space with an unconditional basis \cite[proposition 2.4, p. 69]{Casazza.D.S.Z.2008}. However, this contradicts the fact that 
$\ell_{2}\otimes_{\pi}\ell_{2}$ is not isomorphic to any subspace 
of a Banach space with an unconditional basis \cite[corollary 3.1, p. 58]{Kwapien.Pelski.1970}.\qed
\begin{theorem}
\label{dual.bess.wsc}
Let $\mathcal{G}:=\left( \left( b_{n}^{*},\phi_{n}^{**}\right) \right)
_{n\in \mathbb{N}^{*}}$ be a Besselian frame for $E^{*}$. Then $E^{*}$ is WSC.
\end{theorem}
\textbf{Proof.}
Let $\left( y_{k}^{*}\right)_{k\in\mathbb{N}^{*}}$ be a weakly Cauchy sequence in $E^{*}$. 
Then for each $\phi\in E^{**}$, the limit $\lim_{k\rightarrow \infty}\phi(y_{k}^{*})$ exists. In particular, there exists an element $y_{0}^{*}\in E^{*}$ such that $y_{k}^{*}(x)\rightarrow y_{0}^{*}(x)$
for each $x\in E$. To complete the proof, we show that
$\phi(y_{k}^{*})\rightarrow \phi(y_{0}^{*})$ for each $\phi\in E^{**}$.\\
Let $\phi\in E^{**}$. Since $E^{*}$ is separable, Gelfand theorem (see \cite{Gantmakhter.Smulian.1937}) guarantees the existence of a sequence 
$\left(  x_{r}\right)_{r\in\mathbb{N}^{*}}$ in $E$ such that 
$$\phi(y^{*})=\underset{r\rightarrow \infty}{\lim }y^{*}(x_{r}),
\;\;\;(y^{*}\in E^{*}).$$
For each $n\in\mathbb{N}^{*}$, define  $\alpha_{n}=\underset{k\rightarrow \infty}{\lim }
\phi_{n}^{**}(y_{k}^{*})$. For each $n\in\mathbb{N}^{*}$
and $\beta_{1}, \beta_{2},...,\beta_{n}\in \mathbb{R}$, we have:
\begin{align*}
\left \vert 
 \underset{k=1}{\overset{n}{\sum }}
 \beta_{k}\alpha_{k}
 \right \vert
 &=\underset{i\rightarrow \infty}{\lim }
 \left \vert 
 \underset{k=1}{\overset{n}{\sum }}
 \beta_{k} \phi_{k}^{**}(y_{i}^{*})
 \right \vert
 = \underset{i\rightarrow \infty}{\lim }
 \left \vert 
 \left (
 \underset{k=1}{\overset{n}{\sum }}
 \beta_{k} \phi_{k}^{**}
 \right )(y_{i}^{*})
 \right \vert\\
 &\leq
 \left (\underset{i}{\sup }\Vert y_{i}^{*}\Vert_{E^{*}} \right )
  \left \Vert
 \underset{k=1}{\overset{n}{\sum }}
 \beta_{k} \phi_{k}^{**}
 \right \Vert_{E^{**}}.
 \end{align*}
By Helly's theorem \cite[theorem 4, p. 55]{S.Banach.1932}, there exists 
 $z^{***}\in E^{***}$ such that
 $\left( \alpha_{n}\right)_{n\in\mathbb{N}^{*}}
 =\left(z^{***}(\phi_{n}^{**})\right)_{n\in\mathbb{N}^{*}}$. It follows from lemma \eqref{interchange.E} that
\begin{align*}
 \underset{n=1}{\overset{\infty }{\sum }}
  z^{***}(\phi_{n}^{**})b_{n}^{*}(x_{r}) 
&=
\underset{n=1}{\overset{\infty}{\sum }}
\underset{k\rightarrow \infty}{\lim }
 \phi_{n}^{**}(y_{k}^{*})J_{E}(x_{r})( b_{n}^{*})\\
 &=
 \underset{k\rightarrow \infty}{\lim }
\underset{n=1}{\overset{\infty}{\sum }}
 \phi_{n}^{**}(y_{k}^{*})b_{n}^{*}(x_{r})=y_{0}^{*}(x_{r}),
\end{align*}
then the limit of $\sum_{n=1}^{\infty}
  z^{***}(\phi_{n}^{**})b_{n}^{*}(x_{r}) $ exists as $r$ tends to infinity.
Applying lemmas \eqref{interchange.E} and \eqref{interchange.E**}, we obtain:
\begin{align*}
\underset{k\rightarrow \infty}{\lim }
\phi(y_{k}^{*})
&=
\underset{k\rightarrow \infty}{\lim }
 \underset{r\rightarrow \infty}{\lim }
 y_{k}^{*}(x_{r})
 =
 \underset{k\rightarrow \infty}{\lim }
 \underset{r\rightarrow \infty}{\lim }
 \underset{n=1}{\overset{\infty }{\sum }}
\phi_{n}^{**}(y_{k}^{*}) b_{n}^{*}(x_{r}) \\
&=
\underset{k\rightarrow \infty}{\lim }
\underset{n=1}{\overset{\infty}{\sum }}
\phi_{n}^{**}(y_{k}^{*}) \phi( b_{n}^{*}) 
 = \underset{n=1}{\overset{\infty}{\sum }}
\alpha_{n} \phi( b_{n}^{*})\\
 &= \underset{n=1}{\overset{\infty}{\sum }}
z^{***}(\phi_{n}^{**}) \phi( b_{n}^{*})
=\underset{r\rightarrow \infty}{\lim }
 \underset{n=1}{\overset{\infty }{\sum }}
  z^{***}(\phi_{n}^{**})b_{n}^{*}(x_{r}) \\
 &=\underset{r\rightarrow \infty}{\lim }
 \underset{n=1}{\overset{\infty }{\sum }}
 \underset{k\rightarrow \infty}{\lim }
  \phi_{n}^{**}(y_{k}^{*})J_{E}(x_{r})(b_{n}^{*})
  =\underset{r\rightarrow \infty}{\lim }
 \underset{k\rightarrow \infty}{\lim }
 \underset{n=1}{\overset{\infty }{\sum }}
  \phi_{n}^{**}(y_{k}^{*})b_{n}^{*}(x_{r})  \\
 &=
\underset{r\rightarrow \infty}{\lim }
 \underset{k\rightarrow \infty}{\lim }
 y_{k}^{*}(x_{r})
 =
 \underset{r\rightarrow \infty}{\lim }
 y_{0}^{*}(x_{r})
 =
\phi(y_{0}^{*}).
\end{align*}
This completes the proof.\qed
 \begin{lemma}\cite{kabbaj.karkri.zoubeir.2023}
 \label{E.and.E*.wsc}
Suppose that $E$ has a Besselian frame and both $E$ and $E^{*}$ are WSC.
Then $E$ is reflexive.
 \end{lemma}
 The following corollary extends  the Karlin corollary \cite[p. 984]{S.Karlin.1948} to Banach spaces with Besselian frames.
 \begin{corollary}
 If $E^{*}$ has a Besselian frame and $E^{**}$ is separable, then $E$ is reflexive. 
 \end{corollary}
 \textbf{Proof.}
By theorem \eqref{dual.bess.wsc} and corollary \eqref{bess.dual.bess},
$E^{*}$ and $E^{**}$ are WSC. Lemma \eqref{E.and.E*.wsc} then implies that
 $E^{*}$ is reflexive, and consequently $E$ is reflexive.\qed\\
 
 The following corollary extends  the Karlin corollary \cite[p. 983]{S.Karlin.1948} to Banach spaces with Besselian frames.
\begin{corollary}
Suppose that $E$ has a Besselian frame and $E^{**}$ is separable. Then $E$ is reflexive. 
 \end{corollary}
 \textbf{Proof.} 
By theorem \eqref{bess.dual.WSC} and corollary \eqref{bess.dual.bess}, 
$E^{*}$ is WSC and has a Besselian frame. Since $E^{**}$  is separable,
theorem \eqref{bess.dual.WSC} further implies that $E^{**}$ is WSC. Applying 
 lemma \eqref{E.and.E*.wsc}, we conclude that $E^{*}$ is reflexive.\qed\\
 
It is well-known (see,\cite[proposition 3.5.3, p. 67]{Albiac.Kalton.2006}) that every Banach space with unconditional  basis
has property $(u)$. The following theorem extends this result to Banach spaces 
with Besselian frames (and consequently, to those with unconditional frames).
\begin{theorem}
 \label{property.u}
 Suppose that $\mathcal{F}$ is a Besselian frame for $E$. Then $E$ has the property (u).
 \end{theorem}
 \textbf{Proof.}
 Let $\left( x_{k}\right)_{k\in\mathbb{N}^{*}}$ be a weakly Cauchy sequence in  $E$.
 For each $n\in\mathbb{N}^{*}$, define:
$$\alpha_{n}=\underset{k\rightarrow \infty}{\lim }
b_{n}^{*}(x_{k}).$$
 Since $\mathcal{F}$ is Besselian, then for each  $N,k\in\mathbb{N}^{*}$, we have
$$
\underset{n=1}{\overset{N}{\sum }}\left \vert b_{n}^{*}\left(
x_{k}\right) \right \vert \left \vert y^{*}\left( a_{n}\right) \right \vert
\leq \mathcal{L}_{\mathcal{F}}\left (\underset{k}{\sup }\Vert x_{k}\Vert_{E} \right )\left \Vert y^{*}\right \Vert_{E^{*}},
\;\;\;(y^{*}\in E^{*}).
$$
Taking limits as $k\rightarrow\infty$, we obtain:
$$
\underset{n=1}{\overset{\infty}{\sum }}\left \vert \alpha_{n}\right \vert \left \vert y^{*}\left( a_{n}\right) \right \vert
\leq \mathcal{L}_{\mathcal{F}}\left (\underset{k}{\sup }\Vert x_{k}\Vert_{E} \right )\left \Vert y^{*}\right \Vert_{E^{*}},
\;\;\;(y^{*}\in E^{*}),
$$
which shows the series $\sum_{n}\alpha_{n}a_{n}$ is WUC. For any $y^{*}\in E^{*}$,
lemma \eqref{interchange.E} implies that:
$$
\underset{n=1}{\overset{\infty}{\sum }}
 \alpha_{n}y^{*}(a_{n})
 =
\underset{k\rightarrow \infty}{\lim }
 \underset{n=1}{\overset{\infty }{\sum }}
 b_{n}^{*}(x_{k})y^{*}(a_{n})
 =
 \underset{k\rightarrow \infty}{\lim } 
 y^{*}(x_{k}).
$$
Hence, the sequence $\left( x_{n}-\sum_{i=1}^{n}\alpha_{i}a_{i}\right)_{n\in\mathbb{N}^{*}}$ converges weakly to $0$. That is $E$ has the property (u).\qed
\begin{corollary}
The James space $\mathcal{J}$ has no Besselian frame and consequently, has no unconditional frame.
\end{corollary}
\textbf{Proof.} Since $\mathcal{J}$ fails to have the property (u) \cite[proposition 3.5.4, p. 69]{Albiac.Kalton.2006}, theorem \eqref{property.u} immediately implies that $\mathcal{J}$ has no Besselian frame and consequently, has no unconditional frame.\qed\\

It is well-known (see, \cite[p. 96]{Albiac.Kalton.2006}) that  
 the space $C(\omega^{\omega})$ provides an important example of a Banach space which fails to have
an unconditional basis, despite  its dual space 
having one. The following corollary extends this phenomenon to Banach spaces with Besselian frames.
\begin{corollary}
There exists a Banach space that fails to have a Besselian (and hence unconditional) frame, but its dual space has a Besselian frame.
\end{corollary}
\textbf{Proof.}
Since $C(\omega^{\omega})$ is not isomorphic to $c_{0}$ \cite[p. 284]{I.Singer.I}, it follows from \cite[theorem 4.5.2, p. 96]{Albiac.Kalton.2006}
that $C(\omega^{\omega})$ lacks property (u). By theorem \eqref{property.u}, we conclude that $C(\omega^{\omega})$ cannot have a Besselian frame.\qed  
\begin{corollary}
The  space $C([0,1])$ has no Besselian frame and consequently, has no unconditional frame.
\end{corollary}
\textbf{Proof.}
The space $C([0,1])$ lacks property (u) (see \cite[p. 96]{Albiac.Kalton.2006}). 
By theorem \eqref{property.u}, it follows that $C([0,1])$ cannot have a Besselian frame.\qed
\begin{theorem}\cite[theorem 4.4.7, p. 403]{R.Megginson.1998}
\label{shrinking.basis.equi}
Let $\left( \left( x_{n},x_{n}^{*}\right) \right)
_{n\in \mathbb{N}^{*}}$ be a  basis for  $E$. The following are equivalent:
\begin{enumerate}
\item
The basis $\left( \left( x_{n},x_{n}^{*}\right) \right)
_{n\in \mathbb{N}^{*}}$ is shrinking, i.e., for every $y^{*}\in E^{*}$, 
$$\underset{m\rightarrow \infty}{\lim }
  \Vert y^{*}_{\left \vert span(x_{i}: i>m)\right. }  \Vert_{E^{*}}=0.$$
\item
The sequence $\left( \left(x_{n}^{*}, J_{E}(x_{n})\right) \right)
_{n\in \mathbb{N}^{*}}$ forms a basis for $E^{*}$.
\end{enumerate}
\end{theorem}
\begin{proposition}
\label{shrinking.condition}
Suppose that $\mathcal{F}$ is a Besselian frame for $E$ such that 
$$\underset{m\rightarrow \infty}{\lim }
  \Vert y^{*}_{\left \vert span(a_{i}: i>m)\right. }\Vert_{E^{*}}=0,\;\;(y^{*}\in E^{*}).$$
Then the sequence $\mathcal{F}^{*}$ is a frame for $E^{*}$.
\end{proposition}
\textbf{Proof.}
Assume condition \eqref{shrinking.condition} holds. For any $x\in E$ and $y^{*}\in E^{*}$, we have:
  \begin{align*}
\left \vert \left (y^{* }-\overset{m}{\underset{n=1}{\sum }}
J_{E}(a_{n})\left(y^{*}\right) b_{n}^{*}\right)(x)\right \vert
& =
\left \vert
 y^{*}\left( \overset{\infty }{\underset{n=m+1}{\sum }}
 b_{n}^{* }(x)a_{n}\right) \right \vert
\\
&\leq\Vert y^{*}_{\left \vert span(a_{i}: i>m)\right. }  \Vert_{E^{*}}
\left \Vert \left( \overset{\infty }
{\underset{n=m+1}{\sum }}
b_{n}^{* }(x)a_{n}\right)\right \Vert_{E}
\\
&\leq\mathcal{L}_{\mathcal{F}}
\Vert y^{*}_{\left \vert span(a_{i}: i>m)\right. }  \Vert_{E^{*}}
\Vert x \Vert_{E}.
\end{align*}
Consequently, $\mathcal{F}^{*}$ is a frame for $E^{*}$.\qed\\

Theorem \eqref{shrinking.basis.equi} motivates the following definition.
\begin{definition}
A frame $\mathcal{F}$ for $E$ is called shrinking if its dual sequence $\mathcal{F}^{*}$ forms a frame for $E^{*}$.
\end{definition}
 The following result extends the classical James lemma \cite[p. 520]{James.1950} and \cite[theorem 4.4.21, p. 409]{R.Megginson.1998} to Banach spaces with Besselian frames.
\begin{theorem}
\label{shrinking.copy.l1}
Suppose that  $\mathcal{F}$ is a Besselian frame for $E$. Then $\mathcal{F}$
is shrinking if and only there is no subspace of $E$  isomorphic to 
$\ell_{1}$.  
\end{theorem}
\textbf{Proof.} Suppose that $E$ has a subspace isomorphic to 
$\ell_{1}$, and let $T:\ell_{1}\rightarrow E$ be an isomorphic embedding. Then $T^{*}$ maps $E^{*}$ onto the nonseparable space $\ell_{1}^{*}$ 
\cite[theorem 3.1.22, p. 293]{R.Megginson.1998}, so $E^{*}$ cannot be separable and therefore $\mathcal{F}$ is not shrinking.\\
Assume now that $E$ has no subspace isomorphic to $\ell_{1}$. By \cite[theorem 10, p. 48]{Diestel.1984},
$E^{*}$ has no subspace isomorphic to 
$c_{0}$. For any $y^{*}\in E^{*}$ and $z^{**}\in E^{**}$, the series
 $\sum_{n}\vert z^{**}(y^{*}(a_{n})b_{n}^{*})\vert$ converges. Applying
 \cite[theorem 8, p. 45]{Diestel.1984}, the series $\sum_{n}y^{*}(a_{n})b_{n}^{*}$ is unconditionally convergent in $E^{*}$. This establishes that
 $\mathcal{F}^{*}$ is a frame for $E^{*}$, proving $\mathcal{F}$ is shrinking.\qed\\

The following result extends \cite[corollary 4.4.22, p. 410]{R.Megginson.1998} to Banach spaces with Besselian frames.
\begin{corollary}
For any Banach space $E$, either:
\begin{enumerate}
\item
Every Besselian frame for $E$ is shrinking, or
\item
No Besselian frame for $E$ is shrinking.
\end{enumerate}
 \end{corollary}
\textbf{Proof.} This is an immediate consequence of theorem \eqref{shrinking.copy.l1}.\qed
\begin{theorem}\cite[theorem 6.2, p. 288]{I.Singer.I}
\label{boundedly.basis.equi}
Let $\left( \left( x_{n},x_{n}^{*}\right) \right)
_{n\in \mathbb{N}^{*}}$ be a  basis for  $E$. The following are equivalent:
\begin{enumerate}
\item
The basis is boundedly complete: For any sequence of scalars $(\alpha_{n})_{n\in\mathbb{N}^{*}}$, if 
$$\underset{n\in \mathbb{N}^{*}}{\sup}
 \left \Vert \sum_{i=1}^{n}
 \alpha_{i}x_{i}\right \Vert_{E}<\infty ,$$ 
then the series $\sum_{n}\alpha_{n}x_{n}$ converges in $E$.
\item
For every $\phi^{**}\in E^{**}$, the series 
$\sum_{n}\phi^{**}(x_{n}^{*})x_{n}$ converges in $E$.
\end{enumerate}
\end{theorem}
\begin{proposition}
Suppose that $\mathcal{F}$ is a  Besselian frame for  $E$ satisfying:
For any sequence of scalars $(\alpha_{n})_{n\in\mathbb{N}^{*}}$,
$$
\underset{n\in \mathbb{N}^{*}}{\sup }
 \left \Vert \sum_{i=1}^{n}
 \alpha_{i}a_{i}\right \Vert_{E}<\infty \Longrightarrow \sum_{n}\alpha_{n}a_{n}  \;\;converges.$$
Then, for every $\phi^{**}\in E^{**}$, the series 
$\sum_{n}\phi^{**}(b_{n}^{*})a_{n}$ converges in $E$.
\end{proposition}
\textbf{Proof.}
Let $\phi^{**}\in E^{**}$. Since $\mathcal{F}$ is Besselian, theorem \eqref{F.bess-equi.F*.bess} implies that
$\mathcal{F}^{*}$ is also Besselian. For any $n\in\mathbb{N}^{*}$, we have:
$$\left \Vert \overset{n}{\underset{i=1}{\sum }}
 \phi^{**}(b_{i}^{*})a_{i}\right \Vert_{E}
 \leq 
 \underset{y^{*}\in\mathbb{B}_{E^{*}}}{\sup }
  \overset{n}{\underset{i=1}{\sum }}
 \left \vert\phi^{**}(b_{i}^{*})y^{*}(a_{i})\right \vert
 \leq
 \mathcal{L}_{\mathcal{F}}
\left \Vert \phi^{**} \right \Vert_{E^{**}}. $$
By the hypothesis on $\mathcal{F}$, this implies that the series 
$\sum_{n}\phi^{**}(b_{n}^{*})a_{n}$ converges in $E$.\qed\\

Theorem \eqref{boundedly.basis.equi} motivates the following definition in  \cite{kabbaj.karkri.zoubeir.2023}.
\begin{definition}
A  frame $\mathcal{F}$ for $E$ is called boundedly
complete if for every $x^{**}\in E^{**}$, the series 
$$\sum_{n} x^{**}\left( b_{n}^{*}\right)a_{n}$$
converges in $E$.
\end{definition}
The following result extends the James lemma \cite[p. 520]{James.1950}
 and \cite[theorem 1.c.10, p. 22]{Lindenstraus.Tzafriri.1977} to Banach spaces with Besselian frames.
\begin{theorem}
\label{boundedly.copy.c0}
Suppose that $\mathcal{F}$ is a Besselian frame for $E$. The following are equivalent:
\begin{enumerate}
\item
$\mathcal{F}$ is a boundedly complete  frame for $E$. 
\item
$E$ is weakly sequentially complete.
\item
There is no subspace of $E$ isomorphic to $c_{0}$.
\end{enumerate}
\end{theorem}
\textbf{Proof.}
\begin{itemize}
\item
$(1)\Longrightarrow (2)$: Assume that $\mathcal{F}$ is boundedly complete. Let
 $\left( x_{r}\right)_{r\in\mathbb{N}^{*}}$ be a weakly Cauchy sequence in $E$. For each $n\in\mathbb{N}^{*}$, define:
$$\alpha_{n}=\underset{r\rightarrow \infty}{\lim }
 b_{n}^{*}(x_{r}).$$
Since $\mathcal{F}$ is Besselian, it follows 
 by lemma \eqref{interchange.E} that:
\begin{align}
\label{karlin.1}
 \underset{r\rightarrow \infty}{\lim }
 y^{*}(x_{r})
=
\underset{n=1}{\overset{\infty}{\sum }}
 \alpha_{n}y^{*}\left(a_{n}\right),\;\;\;(y^{*}\in E^{*}).
\end{align}
If the series $\sum_{n}\alpha_{n}a_{n}$ is convergent, then 
$\left( x_{r}\right)_{r\in\mathbb{N}^{*}}$  converges weakly to 
$\sum_{n=1}^{\infty}\alpha_{n}a_{n}$. To verify convergence, observe that for any $n\in\mathbb{N}^{*}$ and $\beta_{1}, \beta_{2},...,\beta_{n}\in \mathbb{R}$:
\begin{align*}
\left \vert 
 \underset{k=1}{\overset{n}{\sum }}
 \beta_{k}\alpha_{k}
 \right \vert
 &=\underset{r\rightarrow \infty}{\lim }
 \left \vert 
 \underset{k=1}{\overset{n}{\sum }}
 \beta_{k} b_{k}^{*}(x_{r})
 \right \vert
 = \underset{r\rightarrow \infty}{\lim }
 \left \vert 
 \left (
 \underset{k=1}{\overset{n}{\sum }}
 \beta_{k} b_{k}^{*}
 \right )(x_{r})
 \right \vert\\
& \leq
 \left (\underset{r}{\sup }\Vert x_{r}\Vert_{E} \right )
  \left \Vert
 \underset{k=1}{\overset{n}{\sum }}
 \beta_{k} b_{k}^{*}
 \right \Vert_{E^{*}}.
 \end{align*}
By Helly's theorem \cite[theorem 4, p. 55]{S.Banach.1932}, there exists 
 $z^{**}\in E^{**}$ such that
 $\left( \alpha_{n}\right)_{n\in\mathbb{N}^{*}}
 =\left(z^{**}( b_{n}^{*})\right)_{n\in\mathbb{N}^{*}}$.
Thus, $\sum_{n}\alpha_{n}a_{n}=\sum_{n}z^{**}( b_{n}^{*})a_{n}$, which converges by the bounded completeness of $\mathcal{F}$. Hence,
 $\left( x_{r}\right)_{r\in\mathbb{N}^{*}}$ converges weakly to
 $\sum_{n=1}^{\infty}\alpha_{n}a_{n}$, proving $E$ is WSC.
\item
$(2)\Longrightarrow (1)$: This implication is clear.
\item
$(3)\Longrightarrow (1)$:
If $E$ has no subspace isomorphic to $c_{0}$, then \cite[theorem.8, p. 45]{Diestel.1984} and theorem \eqref{F.bess-equi.F*.bess} imply that $\mathcal{F}$ is boundedly complete. 
 \item
$(2)\Longrightarrow (3)$: If $E$ is WSC, then 
 by \cite[theorem.8, p. 45]{Diestel.1984}, $E$ has no subspace isomorphic to $c_{0}$.
\end{itemize}\qed
 
The following result extends \cite[corollary 4.4.22, p. 410]{R.Megginson.1998} to Banach spaces with Besselian frames.
\begin{corollary}
For any Banach space $E$, either:
\begin{enumerate}
\item
Every Besselian frame for $E$ is boundedly complete, or
\item
No Besselian frame for $E$ is boundedly complete.
\end{enumerate}
 \end{corollary}
 \textbf{Proof.} This is an immediate consequence of
 theorem \eqref{boundedly.copy.c0}.\qed\\
 
 The following result extends  the James theorem \cite[p. 521]{James.1950}
 and \cite[corollary 4.4.23, p. 410]{R.Megginson.1998} to Banach spaces with Besselian frames.
 \begin{corollary}
 Suppose that $E$ is a Banach space with a Besselian frame. Then the following are equivalent:
\begin{enumerate}
\item
$E$ is reflexive.
\item 
$E$ contains no subspace isomorphic to $\ell_{1}$ or $c_{0}$.
\end{enumerate}
 \end{corollary}
 \textbf{Proof.} This follows immediately from 
 theorems \eqref{shrinking.copy.l1}, \eqref{boundedly.copy.c0}
 and  \cite[theorem 3.10, p. 234]{kabbaj.karkri.zoubeir.2023}.\qed
\section*{Declarations}
\begin{enumerate}
\item 
\textbf{Data Availability Statement.}\\
Not applicable.
\item
\textbf{Competing Interests.}\\
I have nothing to declare.
\item
 \textbf{Funding.}\\
This research did not receive any specific grant from funding agencies in the public, commercial, or not-for-profit sectors.
\item
\textbf{Authors contributions.}\\
The authors equally conceived of the study, participated in its design and coordination, drafted the manuscript, participated in the sequence alignment, 
and read and approved the final manuscript.
\item
\textbf{Ethical approval.}\\
This article does not contain any studies with animals performed by any of the authors.
\end{enumerate}

\bigskip
\bigskip
\end{document}